\documentclass[12pt,leqno]{amsart}
\usepackage[francais]{babel}
\usepackage{amssymb,a4}
\usepackage[all]{xy}
\newcommand{\QQ}{{\mathbb{Q}}}
\newcommand{\ZZ}{{\mathbb{Z}}}

\newcommand\cA{{\mathcal{A}}}

\newcommand\cE{{\mathcal{E}}}
\newcommand\cF{{\mathcal{F}}}
\newcommand\cG{{\mathcal{G}}}
\newcommand\cH{{\mathcal{H}}}
\newcommand\cM{{\mathcal{M}}}
\newcommand\cO{{\mathcal{O}}}

\newcommand\Hom{\operatorname{Hom}}
\newcommand\Irr{\operatorname{Irr}}
\newcommand\Ind{\operatorname{Ind}}

\newcommand\Res{\operatorname{Res}}
\newcommand\Tr{\operatorname{Tr}}
\newcommand\df{\operatorname{def}}

\newcommand\iso{\buildrel \sim\over\longrightarrow}

\newcommand{\ird}{{\sqrt{-3}}}
\newcommand{\zd}{\zeta_3}

\newcommand\fp{{\mathfrak{p}}}
\newcommand\fq{{\mathfrak{q}}}
\newcommand\fn{{\mathfrak{n}}}

\newcommand\bpi{{\boldsymbol\pi}}

\newcommand\cHp{{\cO_\fp\cH}}

\let\eps=\epsilon

\newcommand\CHEVIE{{\sf {Chevie}}}

\newtheorem{thm}{Th\'eor\`eme}[section]
\newtheorem{lem}[thm]{Lemme}
\newtheorem{conj}[thm]{Conjecture}
\newtheorem{cor}[thm]{Corollaire}
\newtheorem{prop}[thm]{Proposition}

\theoremstyle{definition}
\newtheorem{exmp}[thm]{Exemple}
\newtheorem{defn}[thm]{D\'efinition}

\theoremstyle{remark}
\newtheorem{rem}[thm]{Remarque}

\begin{document}

\title{Familles de caract\`eres de groupes de r\'eflexions complexes}

\author{Gunter Malle}
\address{G.M.: Fachbereich Mathematik/Informatik, Universit\"at Kassel, 
Heinrich-Plett-Str.~40, D--34132 Kassel, Germany.}
\makeatletter
\email{malle@mathematik.uni-kassel.de}
\makeatother

\author{Rapha\"el Rouquier}
\address{R.R.: UMR 7586 du CNRS et UFR de Math\'ematiques,
Universit\'e Denis Diderot,
  Case 7012, 2 Place Jussieu, F--75251 Paris Cedex 05, France.}
\makeatletter
\email{rouquier@math.jussieu.fr}
\makeatother


\begin{abstract} 
Nous \'etudions certains types de blocs d'alg\`ebres de Hecke associ\'ees
aux groupes de r\'eflexions complexes qui g\'en\'eralisent les familles
de caract\`eres d\'efinies par Lusztig pour les groupes de Weyl.
Nous d\'eterminons ces blocs pour les groupes de
r\'eflexions spetsiaux et nous \'etablissons un th\'eor\`eme de
compatibilit\'e entre familles et $d$-s\'eries de Harish-Chandra.
\end{abstract}

\maketitle
\pagestyle{myheadings}

\markboth{Malle et Rouquier}{Familles de caract\`eres}

\section{Introduction} \label{sec:intro}

Lusztig a construit une partition
des caract\`eres irr\'eductibles d'un groupe de Coxeter fini $W$ en
{\em familles}, \`a l'aide de la th\'eorie des cellules.
Cette partition appara\^{\i}t naturellement dans le param\'etrage de Lusztig
des caract\`eres unipotents d'un groupe r\'eductif sur un corps fini. \par
Gyoja \cite{Gy} a d\'emontr\'e (cas par cas) que la partition en ``$p$-blocs''
de l'alg\`ebre de Iwahori-Hecke affine la partition en familles et que les deux
partitions co\"{\i}ncident lorsque $p$ est l'unique mauvais nombre premier.
Dans \cite{Ro}, il a ensuite \'et\'e d\'emontr\'e
(par un argument g\'en\'eral utilisant l'alg\`ebre $J$ de Lusztig et une
propri\'et\'e des anneaux bas\'es)
que les blocs de l'alg\`ebre d'Iwahori-Hecke sur un anneau
convenable $\cO$
co\"{\i}ncident avec les familles.
\par
Les groupes de Weyl finis forment une classe particuli\`ere de
groupes de r\'eflexions complexes. Des travaux r\'ecents ont apport\'es
de nombreux indices que les groupes de r\'eflexions complexes (au moins
certains d'entre eux) jouissent de propri\'et\'es analogues aux groupes de
Weyl (voir par exemple \cite{spets,MaS}). On dispose en particulier
d'alg\`ebres du type Iwahori-Hecke.
Pour le moment, il
n'existe pas de d\'efinition de cellules pour les groupes de
r\'eflexions complexes. On ne peut donc utiliser
l'approche de Lusztig pour d\'efinir les
familles de caract\`eres, mais celle de \cite{Ro} d\'ecrite plus haut
est differente~: on d\'efinit
les familles comme blocs de l'alg\`ebre de Iwahori-Hecke sur
un certain anneau $\cO$.\par

La premi\`ere partie de cet article consiste en la d\'etermination, cas par
cas, des $\cO$-blocs des alg\`ebres de Hecke
pour les groupes de r\'eflexions exceptionnels spetsiaux, c'est \`a
dire, pour les~18 groupes de r\'eflexions irr\'eductibles
$$G_4,\ G_6,\ G_8,\ G_{14},\ G_{23,\ldots,30},\ G_{32,\ldots,37}$$
dans la notation de Shephard-Todd (ceci couvre en particulier le
cas du groupe de Coxeter de type $H_4$ qui ne pouvait \^etre trait\'e par
les m\'ethodes de \cite{Ro}). Le cas des s\'eries infinies a \'et\'e
trait\'e par Brou\'e et Kim \cite{BrKi}. Nous calculons aussi (partiellement)
les matrices de d\'ecomposition pour chaque nombre premier mauvais.
\par

Dans un second temps, nous d\'emontrons
(th\'eor\`eme~\ref{fam=bloc}) un r\'esultat concernant les caract\`eres
unipotents de groupes r\'eductifs finis~:
\`a l'int\'erieur d'une $d$-s\'erie de Harish-Chandra, les familles
d\'efinies par l'alg\`ebre de Hecke de la s\'erie sont les intersections
des familles de Lusztig avec la $d$-s\'erie concern\'ee.
Ce r\'esultat (\'etabli cas par cas)
peut \^etre vu comme une g\'en\'eralisation
des r\'esultats de \cite{Ro}, qui ne concernent que le cas
$d=1$ et la s\'erie principale.

Dans le cas particulier des $1$-s\'eries de Harish-Chandra,
le th\'eor\`eme~\ref{fam=bloc} appara\^{\i}t comme une variante
de la conjecture de Lusztig \cite[(25.3)]{Lu3}
sur les cellules bilat\`eres dans les alg\`ebres de Hecke usuelles
avec des param\`etres in\'egaux.

\par
Les calculs ont \'et\'e effectu\'es \`a l'aide du logiciel
\CHEVIE \cite{chev}, en utilisant en particulier les fonctions sur les
groupes de r\'eflexions complexes impl\'ement\'ees par Jean Michel. Nous
lui adressons nos remerciements pour avoir fourni cet outil essentiel.
Nous remercions aussi J\"urgen M\"uller pour ses
calculs de d\'ecomposition de repr\'esenta\-tions (en particulier,
pour $G_{30}=W(H_4)$) avec les logiciels {\sf MeatAxe} et
{\sf VectorEnumerator}.

\section{Alg\`ebres de Hecke de groupes de r\'eflexions complexes et familles}
\label{sec:theorie}

Nous commen\c{c}ons par rappeler la d\'efinition des alg\`ebres de Hecke
associ\'ees aux groupes de r\'eflexions complexes et
certaines de leurs propri\'et\'es. Ensuite, nous donnons quelques r\'esultats
qui seront utilis\'es plus tard pour la d\'etermination des $\cO$-blocs.

\subsection{Alg\`ebres de Hecke}
\subsubsection{}
Soit $k$ un corps de nombres et $W$ un groupe de r\'eflexions complexe
fini sur le $k$-espace vectoriel $V$.
Soit $\cH=\cH(W,x)$ l'alg\`ebre de
Hecke de $W$ sur l'anneau des polyn\^omes de Laurent
$A=\ZZ[x,x^{-1}]$.

Rappelons sa construction par Brou\'e--Malle--Rouquier \cite{BMR}.
Soit $\cA$ l'ensemble
des hyperplans de r\'eflexion de $W$ et $B_W$ le groupe de tresses
de $W$, {\em i.e.}, le groupe fondamental de
$$X=(V-\bigcup_{H\in\cA}H)/W$$
relatif \`a un point base $x_0$. Alors, $\cH$ est le quotient
de l'alg\`ebre de groupe $A[B_W]$ par l'id\'eal
engendr\'e par $(\sigma_H-x)(1+\sigma_H+\cdots+\sigma_H^{d-1})$, o\`u
$H$ d\'ecrit $\cA$,
$d$ est l'ordre du fixateur de $H$ dans $W$,
$s$ est la r\'eflexion autour de $H$ de d\'eterminant $\exp(2i\pi/d)$ et
$\sigma_H$ est un $s$-g\'en\'erateur de la monodromie autour de $H$.

La structure
de ces alg\`ebres de Hecke reste encore myst\'erieuse. N\'eanmoins
on conjecture qu'elles partagent beaucoup des propri\'et\'es des alg\`ebres
d'Iwahori-Hecke usuelles. 

\subsubsection{}
Soit $W'$ un sous-groupe parabolique de $W$. Dans \cite[\S 4.C]{BMR}, on
construit une injection de $\cH(W')$ dans $\cH$. Cette 
injection est d\'efinie \`a automorphisme int\'erieur de $\cH$ pr\`es.

\begin{conj} \label{conj1}
 Pour tout sous-groupe parabolique $W'$ de $W$, l'alg\`ebre $\cH$ est un
 module libre de rang $[W:W']$ sur $\cH(W')$.
\end{conj}

Notons que cette conjecture est ind\'ependante du choix
du plongement.
On conjecture en particulier que $\cH$ est un module libre de 
rang $|W|$ sur $A$.

On notera $A'\cH$ l'alg\`ebre $A'\otimes_A\cH$, pour $A'$ une $A$-alg\`ebre
commutative.

Supposons dor\'enavant la conjecture \ref{conj1} correcte. Alors,
$\cH$ est une d\'eformation de $\ZZ[W]$~: on a un isomorphisme
\begin{equation}  \label{bij}
  \cH\otimes_A \left(A/(x-1)\right)\iso \ZZ[W]\,
\end{equation}
(ce morphisme est celui par lequel se factorise l'application canonique
$\ZZ[B_W]\to \ZZ[W]$ provenant du rev\^etement galoisien
$V-\bigcup_{H\in\cA}H\to X$).

Soit $K$ une extension galoisienne finie de $k(x)$ telle que
$K\cH$ et $K[W]$ sont semi-simples d\'eploy\'ees (d'apr\`es Benard et Bessis
\cite{Bes}, l'alg\`ebre $k[W]$ est d\'eploy\'ee et d'apr\`es Malle
\cite[Thm.~5.2]{MaR},
l'alg\`ebre $\cH\otimes_A k(x^{1/\mu(k)})$ est d\'eploy\'ee, o\`u $\mu(k)$
est l'ordre du groupe des racines de l'unit\'e de $k$).

Soit $B$ la cl\^oture int\'egrale de $A$ dans $K$ et
$\fn$ un id\'eal premier de $B$ au-dessus de $(x-1)$.
On d\'eduit de (\ref{bij}) un isomorphisme
$(B_\fn/\fn)\cH\iso (B_\fn/\fn)[W]$. Puisque 
$(B_\fn/\fn)\cH$ est semi-simple et d\'eploy\'ee, l'application
de d\'ecomposition est une bijection entre l'ensemble $\Irr(\cH)$ des
caract\`eres irr\'eductibles de $K\cH$ et l'ensemble des caract\`eres
irr\'eductibles de $(B_\fn/\fn)\cH$. De m\^eme, on a une bijection
entre l'ensemble $\Irr(W)$ des caract\`eres irr\'eductibles de $K[W]$
et l'ensemble des caract\`eres irr\'eductibles de $(B_\fn/\fn)[W]$.

Par cons\'equent, on obtient une bijection
\begin{equation} \label{carbij}
  \Irr(\cH)\iso\Irr(W)\,.
\end{equation}
Dans la suite, nous identifierons les 
caract\`eres irr\'eductibles de $W$ et de $\cH$ via cette bijection. 

Cette bijection d\'epend du choix de l'id\'eal premier $\fn$.
Les id\'eaux premiers au dessus de $(x-1)$ sont permut\'es par le groupe
de Galois de $K/k(x)$. Soit maintenant $R$ une sous-$A$-alg\`ebre de $K$
invariante par le groupe de Galois. La partition des caract\`eres
irr\'eductibles de $K\cH$ selon les blocs de $R\cH$ est invariante par
action du groupe de Galois~: par cons\'equent, la partition correspondante
de $\Irr(W)$ est ind\'ependante du choix de $\fn$.
\par

\subsubsection{}
Dans Brou\'e--Malle--Michel \cite[Thm.~2.1]{spets}, il est d\'emontr\'e qu'il
existe au plus une
forme sym\'etrisante sur $\cH$ v\'erifiant certaines propri\'et\'es~:
nous l'appellerons {\em forme sym\'etrisante canonique}.

\begin{conj} \label{conj2}
 La forme sym\'etrisante canonique sur
 $\cH(W)$ existe et sa restriction \`a toute sous-alg\`ebre parabolique
 $\cH(W')$ est la forme sym\'etrisante canonique.
\end{conj}

Nous supposons aussi que la forme sym\'etrisante construite dans \cite{MaD}
pour les groupes de r\'eflexions complexes de rang $2$ et
dans \cite{MaG} pour des groupes de rang sup\'erieur est la forme canonique.

Rappelons qu'on ne sait pas, en g\'en\'eral, que $\cH$ est une alg\`ebre
sym\'etrique.

\smallskip

Soit $t:\cH\rightarrow A$ une forme sym\'etrisante v\'erifiant la
conjecture \ref{conj2}.
Alors, on a une d\'ecomposition de l'extension de $t$ \`a $K$~:
$$t = \sum_{\chi\in\Irr(W)} \frac{1}{c_\chi} \chi$$
o\`u nous avons identifi\'e les caract\`eres de $K\cH$ et de $W$.
Les \'el\'ements de Schur $c_\chi$ sont dans $B$ (cf. Geck--Rouquier
\cite[Prop.~4.2]{GeRo}).

\subsubsection{}
La conjecture suivante est vraie lorsque $W$ est un groupe de
Weyl \cite[Thm. 5.2]{GeRo}. Le bon comportement du centre
par extension des scalaires est une propri\'et\'e importante pour l'\'etude
des repr\'esentations d'une alg\`ebre sym\'etrique, car elle est \'equivalente
\`a la propri\'et\'e correspondante pour les fonctions de classe sur
l'alg\`ebre.

\begin{conj} \label{conj3}
 Soit $R$ une $A$-alg\`ebre commutative locale de corps r\'esiduel
 $F$. Alors, le morphisme canonique $Z(R\cH)\to Z(F\cH)$ est surjectif.
\end{conj}

\subsection{D\'etermination des $\cO$-blocs}

Dans la suite, nous supposerons que $W$ v\'erifie les conjectures \ref{conj1},
\ref{conj2} et \ref{conj3}.

Les conjectures~\ref{conj1} et~\ref{conj2} ont \'et\'e v\'erifi\'ees
pour tous les groupes de r\'eflexion complexes imprimitifs ainsi que
pour des groupes primitifs d'ordre petit, via des m\'ethodes
informatiques.
Quant \`a la conjecture~\ref{conj3}, elle n'a pas \'et\'e, \`a notre
connaissance, v\'erifi\'ee pour des groupes de r\'eflexion non
r\'eels.\par

\subsubsection{}
Soit $\cO=B[(1+xP(x))^{-1}]_{P\in\ZZ[x]}$.

Nous identifierons comme \`a l'accoutum\'ee les blocs d'une alg\`ebre $C$
sur un anneau commutatif
int\`egre $R$ avec l'ensemble des caract\`eres irr\'eductibles de
$F\otimes_R C$, o\`u $F$ est le corps des fractions de $R$.

\begin{defn} \label{familles}
Les {\em familles de $\Irr(W)$} sont les blocs de $\cO\cH(W)$
(via l'identification~(\ref{carbij})).
\end{defn}

D'apr\`es \cite{Ro}, cette d\'efinition co\"{\i}ncide avec
la d\'efinition de Lusztig en termes de combinatoire de Kazhdan-Lusztig,
pour les groupes de Weyl et plus g\'en\'eralement
pour les groupes de Coxeter finis n'ayant pas de facteur de type $H_4$
(une propri\'et\'e de positivit\'e des constantes de structure de
la base de Kazhdan-Lusztig est requise).
Nos r\'esultats montrent (cf. corollaire~\ref{family=block}) que c'est en fait
vrai pour tous les groupes de Coxeter finis.

\smallskip
On d\'efinit de mani\`ere inductive la notion de {\em caract\`ere constructible}
de $W$~: le caract\`ere trivial du groupe trivial est constructible.
Les caract\`eres constructibles de $W$ sont les \'el\'ements
de la famille g\'en\'eratrice minimale du sous-mono\"{\i}de de $\ZZ\Irr(W)$
engendr\'e par les $p_\cF(\Ind_{W'}^W \phi)$, o\`u $W'$ est un sous-groupe
parabolique propre de $W$, $\phi$ un caract\`ere constructible de $W'$,
$\cF$ une famille de $\Irr(W)$ et $p_\cF$ la projection sur $\ZZ\cF$.
Les caract\`eres constructibles sont donc des caract\`eres de 
$\cO\cH$-modules projectifs.

\smallskip
Rappelons la construction des caract\`eres constructibles au sens de Lusztig.
C'est la donn\'ee, pour chaque groupe de Coxeter fini, d'un ensemble
minimal de caract\`eres v\'erifiant les propri\'et\'es suivantes~:

\begin{itemize}
\item le caract\`ere trivial est constructible
\item pour $\phi$ un caract\`ere constructible d'un sous-groupe parabolique
$W'$ de $W$, alors le $J$-induit de $\phi$ est un caract\`ere constructible
de $W$
\item le produit tensoriel par le caract\`ere signe $\eps$ d'un
caract\`ere constructible est encore un caract\`ere constructible.
\end{itemize}

\begin{prop}
 Lorsque $W$ est un groupe de Coxeter fini, un caract\`ere
 est constructible (au sens d\'efini ci-dessus) si et seulement si
 il est constructible au sens de Lusztig.
\end{prop}

\begin{proof}[Preuve]
Le caract\`ere trivial est seul dans sa famille et
le caract\`ere signe $\eps$ est constructible.
La $J$-induction est une induction suivie par une projection
sur une union de familles. Par
cons\'equent, un caract\`ere constructible au sens de Lusztig est
une somme de caract\`eres constructibles.

Il suffit de d\'emontrer la proposition pour $W$ irr\'eductible, ce que
nous supposons maintenant.

Si $W$ n'est pas de type $F_4$ ou $E_n$, alors les caract\`eres
constructibles au sens de Lusztig forment une base du sous-mono\"{\i}de
des caract\`eres form\'e des caract\`eres des $\cO\cH$-modules
projectifs \cite[Thm. 2]{Ro}. Or, les
caract\`eres constructibles sont des caract\`eres de $\cO\cH$-modules
projectifs. Par cons\'equent, les caract\`eres constructibles sont
les caract\`eres constructibles au sens de Lusztig.

Le cas o\`u $W$ est de type $F_4$ ou $E_n$ se v\'erifie directement
par calcul de tous les induits de caract\`eres constructibles de
sous-groupes paraboliques stricts maximaux.
\end{proof}

\subsubsection{}
Soit $p$ un nombre premier,
$\fp$ un id\'eal premier de $\cO$ au-dessus de $p$ et $k_\fp$ le corps
des fractions de $\cO_\fp/\fp$.\par

Nous renvoyons le lecteur \`a Geck--Rouquier \cite[\S 2]{GeRo} pour le
traitement des applications de d\'ecomposition (extension au cas d'un anneau
de base int\'egralement clos des r\'esultats classiques pour un anneau
de valuation discr\`ete).

On dispose d'une application de d\'ecomposition du groupe de Grothendieck
$G_0(K\cH)\!$ de la cat\'egorie des $K\cH$-modules de type fini
vers $G_0(k_\fp\cH)$ et
d'une application de Cartan $c$ du groupe de Grothendieck $K_0(k_\fp\cH)$ de
la cat\'egorie des $k_\fp\cH$-modules projectifs de type fini
vers $G_0(K\cH)$ qui est
duale de l'application de d\'ecomposition. L'application de
Cartan est injective (on utilise ici la conjecture 2.3) \cite[Prop.~3.1]{GeRo}.
Nous appellerons caract\`ere
virtuellement projectif (de $\cHp$) un \'el\'ement de $G_0(K\cH)$
dans l'image de l'application de Cartan.
Lorsque cet \'el\'ement provient d'un $k_\fp\cH$-module projectif, nous
dirons que c'est un caract\`ere projectif.

Les familles de caract\`eres de $W$ sont les sous-ensembles minimaux de
$\Irr(W)$ qui sont r\'eunion de blocs de $\cHp$ pour tout nombre premier
$p$.

\smallskip
Rappelons quelques propri\'et\'es classiques utilisant
la structure sym\'etrique de $\cHp$.

\begin{lem} \label{props}
 \begin{itemize}
 \item[(a)] Soit $\phi$ un caract\`ere virtuellement
 projectif de $\cHp$. Alors,
 $$\sum_{\chi\in\Irr(W)} \frac{\langle \phi,\chi\rangle}{c_\chi}\in\cO_\fp\,.$$
 \item[(b)] $\chi$ est seul dans son bloc $\iff$  $c_\chi\notin\fp$ $\iff$
 $\chi$ est un caract\`ere projectif.
\par
 \item[(c)] Si deux caract\`eres irr\'eductibles $\chi$ et $\psi$ sont dans le
 m\^eme bloc de $\cHp$, alors ils sont dans le m\^eme bloc de $\cO_\fp[W]$.
 Tout caract\`ere projectif de $\cO_\fp[W]$ est un caract\`ere projectif
 de $\cHp$. \par
 \item[(d)] Soit $W'$ un sous-groupe parabolique de $W$. Alors,
 la restriction d'un caract\`ere projectif de $\cHp$ \`a 
 $\cO_\fp\cH(W')$ est un caract\`ere projectif. De m\^eme,
 l'induit d'un caract\`ere projectif de $\cO_\fp\cH(W')$ est un caract\`ere
 projectif de $\cHp$.
 \end{itemize}
\end{lem}

\begin{proof}[Preuve]
 L'assertion (a) est \cite[Prop.~4.4]{GeRo}.

 \smallskip
 Si $\chi$ est un caract\`ere projectif, il r\'esulte de (a) que
 $c_\chi\not\in\fp$.

 Soit $C$ le bloc de $\cO_\fp\cH$ contenant $\chi$.
 Si $c_\chi$ est inversible dans $\cO_\fp$, alors l'idempotent central
 $e_\chi$ associ\'e \`a $\chi$ est dans $\cHp$, donc $\chi$ est seul dans
 son bloc, $k_\fp C$ est une alg\`ebre simple et
 $\chi$ est un caract\`ere projectif.

 Si $\chi$ est seul dans son bloc, alors $Z(C)$ est de rang $1$
 sur $\cO_\fp$, donc $Z(k_\fp C)$ est de dimension $1$
 (conjecture~\ref{conj3}). Par cons\'equent, le socle de $k_\fp C$ est nul,
 donc $k_\fp C$ est une alg\`ebre simple et $\chi$ est un caract\`ere
 projectif.

 Si $\chi$ est un caract\`ere projectif, alors $\chi$ est seul dans son bloc.
 Ceci ach\`eve la preuve de (b).

 \smallskip
 Soit $\fq$ un id\'eal premier de $B$ contenant $\fp$
 et $\fn$ ($\fn$ est un id\'eal premier de $B$ au-dessus de $(x-1)$).
 La surjection canonique
 $$(B_\fq)\cH\twoheadrightarrow (B_\fq/\fn)\cH\iso (B_\fq/\fn)[W]$$
 envoie un bloc sur un bloc.
 Puisqu'un bloc de $(B_\fq)\cH$ est une somme de
 blocs de $(B_\fp)\cH=\cO_\fp\cH$ et que les blocs de
 $(B_\fq/\fn)[W]$ co\"{\i}ncident
 avec ceux de $\cO_\fp[W]$, on obtient la premi\`ere partie de (c).
 La seconde partie de (c) r\'esulte de la transitivit\'e des applications
 de d\'ecomposition \cite[Prop.~2.12]{GeRo}.

 \smallskip
 Puisque $\cHp$ est un $\cO_\fp\cH(W')$-module projectif par
 conjecture~\ref{conj1}, les foncteurs
 d'induction et de restriction entre ces alg\`ebres sont exacts, d'o\`u (d).
\end{proof}

\smallskip
Soit $W'$ un sous-groupe parabolique de $W$. On consid\`ere le morphisme
canonique
$$\cH(W)\to \cH(W)\otimes_{\cH(W')}
\cH(W)\to \cH(W)$$
composition de l'unit\'e et de la counit\'e provenant de
la paire biadjointe
$(\Res_{\cH(W')}^{\cH(W)},\Ind_{\cH(W')}^{\cH(W)})$.
On note $\Tr_{\cH(W')}^{\cH(W)}1\in Z(\cH)$ l'image de $1$ par
ce morphisme. Alors, on a le lemme classique suivant~:

\begin{lem}
\label{projrela}
L'action de 
$\Tr_{\cH(W')}^{\cH(W)}1$ sur un $K\cH$-module simple de caract\`ere
$\chi$ est donn\'ee par l'\'el\'ement de $\cO$
$$c_{\chi}\sum_{\psi\in\Irr(W')}
\frac{\langle \Res_{W'}^W\chi,\psi\rangle}{c_{\psi}}.$$

Soit $B$ un bloc de $\cO_\fp\cH(W)$. Alors, l'image dans $k_\fp$ du scalaire
par lequel $\Tr_{\cH(W')}^{\cH(W)}1$ agit sur un $KB$-module simple
est ind\'ependante du module simple. En outre,
$B$ est projectif relativement \`a $\cO_\fp\cH(W')$ si et seulement si
ce scalaire est inversible.
\end{lem}

\begin{proof}[Preuve]
Pour la formule, voir par exemple \cite[Lemma 9.4.6]{GP}.

Pour le reste, on utilise la propri\'et\'e de
$Z(k_\fp B)$ d'\^etre local.
Un \'el\'ement de $Z(B)$ est inversible si et seulement si
son image dans le corps r\'esiduel de $Z(k_\fp B)$ est non nulle, ce qui
explique la derni\`ere assertion.
\end{proof}

\smallskip
Soit $S(V)$ l'alg\`ebre sym\'etrique de l'espace vectoriel $V$. Alors $S(V)$ est
un $k[W]$-module gradu\'e. L'alg\`ebre des coinvariants de $W$ est le
quotient de
$S(V)$ par l'id\'eal engendr\'e par les invariants de degr\'e positif,
$R=S(V)/(S(V)_+^W)$. C'est encore un $k[W]$-module gradu\'e. Le polyn\^ome de
Poincar\'e $P(W)\in\ZZ[x]$ de $W$ est la dimension gradu\'ee de $R$. Pour
$\chi\in\Irr(W)$ le {\em degr\'e fant\^ome} $R_\chi\in\ZZ[x]$ est
d\'efini comme la multiplicit\'e gradu\'ee de $\chi$ dans le caract\`ere de $R$.
Notons $N(\chi)=\frac{dR_\chi}{dx}|_{x=1}$. Pour $\chi\in\Irr(W)$, nous notons
$\chi^*$ le caract\`ere conjugu\'e.

Le r\'esultat suivant est utile pour la d\'etermination des blocs~:

\begin{lem} \label{central}
 L'application
 $$\Irr(W)\longrightarrow\ZZ,\qquad \chi\mapsto (N(\chi)+N(\chi^*))/\chi(1)\,,$$
 est constante sur les blocs de $\cO\cH$.
\end{lem}

\begin{proof}[Preuve]
D'apr\`es Brou\'e--Michel \cite[Prop.~4.18]{BrMi}, on a
$$\chi(T_\bpi)=\chi(1)x^{N+N^*-(N(\chi)+N(\chi^*))/\chi(1)}\,,$$
o\`u $\bpi$ est la classe du lacet $t\mapsto \exp(2i\pi t)x_0$ (un 
\'el\'ement du centre de $B_W$), $T_\bpi$ son image dans $\cH$ (voir aussi
Brou\'e--Malle--Michel \cite[Prop.~6.7(3)]{spets}), $N$ est le nombre
d'hyperplans de r\'eflexion de $W$
et $N^*$ le nombre de r\'eflexions. On d\'eduit alors le lemme
de la propri\'et\'e de la restriction de $\frac{1}{\chi(1)}\chi$
modulo $\fp$ au centre d'un bloc d'\^etre ind\'ependante du caract\`ere
irr\'eductible $\chi$ du bloc.
\end{proof}


\subsection{Groupes spetsiaux}

Suivant \cite{MaG} nous dirons que le groupe de r\'eflexions complexe $W$ est
{\em spetsial} si tous les \'el\'ements de Schur $c_\chi$ de $\cH(W)$ sont
dans $A$. Les
groupes de r\'eflexions spetsiaux sont d\'etermin\'es dans \cite[8A]{MaG}.
Notre but est la d\'etermination des familles de $W$ et des matrices de
d\'ecomposition de
$\cO_\fp\cH$. D'apr\`es le lemme~\ref{props}(b),
l'application de d\'ecomposition
est l'identit\'e si $c_\chi\notin \fp$ pour tout $\chi$.
Nous dirons que $p$ est
un {\em nombre premier mauvais} pour $W$ si $p$ divise un des $c_\chi$,
pour un $\chi\in\Irr(W)$. Il nous suffit donc de consid\'erer les matrices de
d\'ecomposition pour des id\'eaux $\fp$ au-dessus des mauvais nombres
premiers. \par  

Dans la table~\ref{mauvais} nous indiquons les mauvais nombres premiers pour
les groupes de r\'eflexions spetsiaux exceptionnels (i.e., pour les groupes
qui ne font pas partie des s\'eries infinies).

\begin{table}[htbp] \caption{Mauvais nombres premiers} \label{mauvais}
\[\begin{array}{|l|c|c|l|c|}
\hline
 \hfil W& \text{mauvais }p& & \hfil W& \text{mauvais }p\cr
\hline
 G_4& 2,3&          & G_{28}=W(F_4)& 2,3\cr
 G_6& 2,3&          & G_{29}& 2,5\cr
 G_8& 2,3&          & G_{30}=W(H_4)& 2,3,5\cr
 G_{14}& 2,3&       & G_{32}& 2,3,5\cr
 G_{23}=W(H_3)& 2,5&& G_{33}& 2,3\cr
 G_{24}& 2,7&       & G_{34}& 2,3,7\cr
 G_{25}& 2,3&       & G_{35}=W(E_6)& 2,3\cr
 G_{26}& 2,3&       & G_{36}=W(E_7)& 2,3\cr
 G_{27}& 2,3,5&     & G_{37}=W(E_8)& 2,3,5\cr
\hline
\end{array}\]
\end{table}

\section{Description des familles et des matrices de d\'ecomposition}
\label{sec:explicit}

\subsection{D\'emarche algorithmique}\label{subsec:algo}

Nous allons donner les matrices de d\'ecomposition de $\cO_\fp\cH(W)$ pour
tous les groupes de r\'eflexions complexes spetsiaux irr\'eductibles
$W$, \`a des ind\'eterminations
pr\`es pour $W=G_{29}$, $p=2$, $W=G_{32}$, $p=2,3$ et $W=G_{34}$, $p=3$.
Le cas des groupes de Weyl $F_4$, $E_6$ et $E_7$ a \'et\'e r\'esolu
par Gyoja \cite{Gy}.

La notation utilis\'ee pour les caract\`eres co\"{\i}ncide avec celle
du syst\`eme \CHEVIE \cite{chev}.

\smallskip
Commen\c cons par d\'ecrire
notre d\'emarche algorithmique.

Nous proc\'edons par induction.
Soit $W$ un groupe de r\'eflexions
complexe et $p$ un nombre premier mauvais. Les $p$-blocs de $W$ se
d\'eterminent imm\'ediatement \`a partir de la table de caract\`eres de $W$.
Supposons que
les degr\'es fant\^omes $R_\chi$ et les \'el\'ements de Schur $c_\chi$ sont
connus ainsi que les matrices de d\'ecomposition 
correspondantes pour les alg\`ebres de Hecke de tous
les sous-groupes paraboliques stricts (maximaux) de $W$.
Alors une approximation de la matrice de d\'ecomposition est construite par
la m\'ethode suivante~:

\begin{itemize}
\item[(1)] Une partition de $\Irr(W)$ en union de blocs est donn\'ee
par l'intersection des $p$-blocs de $W$ avec les sous-ensembles sur lesquels
$(N(\chi)+N(\chi^*))/\chi(1)$ est constant (cf. lemmes~\ref{central} et
\ref{props}(c)). De plus, chaque $\chi$ avec $c_\chi\notin \fp$ forme
un bloc en lui-m\^eme (cf. lemme~\ref{props}(b)).

\item[(2)] Une approximation des caract\`eres projectifs ind\'ecomposables
s'obtient comme suit. Soit $P$ l'ensemble des induits des
caract\`eres projectifs ind\'ecomposables des sous-alg\`ebres
paraboliques maximales, coup\'es par les approximations des blocs obtenus \`a
la premi\`ere \'etape. On construit alors la famille $U$ g\'en\'eratrice du
sous-mono\"{\i}de de $\ZZ\Irr(W)$ engendr\'e par $P$.

\item[(3)] On d\'efinit une relation d'\'equivalence comme cl\^oture transitive
de la relation suivante sur $\Irr(W)$:
deux caract\`eres $\chi_1,\chi_2$ sont reli\'es si ils
sont dans une m\^eme partie de la partition construite en (1) et si de plus il
existe un caract\`ere $\phi\in U$ tel que $\chi_1$ et $\chi_2$
interviennent dans $\phi$. On obtient ainsi un raffinement de la partition
construite en (1) o\`u les parties sont toujours des unions de blocs.

\item[(4)] Par construction, les \'el\'ements de $U$ sont des caract\`eres
projectifs. Supposons que $\phi\in U$ n'est pas ind\'ecomposable.
Alors, il existe des caract\`eres non nuls $\phi_1$ et $\phi_2$ avec
$\phi_1+\phi_2=\phi$ et $\sum_\chi \langle \phi_1,\chi\rangle/c_\chi\in\cO_\fp$
(cf. lemme~\ref{props}(a)). Par cons\'equent,
si la somme pr\'ecedente n'est contenue dans $\cO_\fp$ pour aucun
sous-caract\`ere propre de $\phi$, alors $\phi$ est ind\'ecomposable.
\end{itemize}

Les points~(1)--(4) suffisent pour obtenir les matrices de
d\'ecomposition pour tous les
mauvais nombres premiers lorsque $W=G_i$ avec $i\in\{4,6,23,26,27,33\}$. Pour
les autres groupes, nous aurons besoin d'arguments suppl\'ementaires.

Les \'etapes d\'ecrites ci-dessus peuvent toutes \^etre effectu\'ees
manuellement. N\'eanmoins, nous avons utilis\'e la librairie \CHEVIE\
pour organiser et manipuler les tables plus ais\'ement et pour r\'eduire
les risques d'erreurs. Par exemple, pour $W=G_{24}$, les instructions
suivantes ont \'et\'e utilis\'ees:

\begin{verbatim}
W:=ComplexReflectionGroup(24);            ## creation de W
ti:=CharTable(W).irreducibles;;           ## recupere Irr(W)
Aa:=[0,42,28,14,28,14,14,28,24,18,21,21]; ## liste des val. de (a+A)
W1:=ReflectionSubgroup(W,[1,2]);          ## cree un ssgp de reflexion
ind1:=InductionTable(W1,W);               ## table d'ind. pour W1<W
for j in [0..100] do if Position(Aa,j)<>false then
  Display(ind1,rec(charsGroup:=Filtered([1..Length(ti)],i->Aa[i]=j)));
fi; od;
##  affiche la table d'induction, divisee suivant la valeur de a+A
\end{verbatim}

Pour utiliser l'information sur les $p$-blocs de $W$, les instructions
GAP suivantes
\begin{verbatim}
t:=CharTable(W);                ##  recupere la table de car. de W
pr:=PrimeBlocks(t,p).block;     ##  determine les p-blocs
\end{verbatim}
ont \'et\'e utilis\'ees.

\goodbreak

\subsection{Types de familles}

De nombreux blocs de $\cO\cH(W)$ pr\'esentent des matrices de d\'ecomposition
identiques pour des groupes $W$ diff\'erents.
Aussi, nous introduisons ici les types rencontr\'es plus d'une fois.
Nous indiquons la matrice de d\'ecomposition pour
chaque nombre premier $p$ o\`u elle n'est pas l'identit\'e.
Nous indiquons aussi le plus petit coefficient non nul $f_\chi$ de $c_\chi$
et dans la colonne $c$, les caract\`eres constructibles. 

\smallskip

\hbox{\hskip 1cm
\vbox{
\vbox{\offinterlineskip\halign{\vrule height11pt depth2.5pt width0pt
         $#$\hfil\quad\vrule\  & $#$\hfil\quad\vrule\ &
         $#$\hfil\quad\vrule\  & $#$\hfil\quad\vrule \cr
    & \hfil f_\chi\hfil & \hfil p\hfil & \hfil c \hfil\cr
\noalign{\hrule}
\chi_1 & \frac{1}{2}p^n(1+\sqrt{1-4p^{-n}}) & 1 & 1\cr
\chi_2 & \frac{1}{2}p^n(1-\sqrt{1-4p^{-n}}) & 1 & 1\cr
\noalign{\hrule}
      }}
\vskip 1pc
Type $\alpha_{1^2}(p^n)$}
\hskip -7cm
\vbox{
 \vbox{\offinterlineskip\halign{\vrule height10pt depth2.5pt width0pt
         $#$\hfil\quad \vrule\  & \hfil$#$\ \vrule\ &
         $#$\hfil\quad\vrule\  & $#$\hfil\quad\vrule \cr
          & \hfil f_\chi\hfil & \hfil p=2\hfil & \hfil c\hfil\cr
\noalign{\hrule}
\chi_1 & 2(1-\sqrt{-1})& 1 & 1\cr
\chi_2 &              2& 1 & 1\cr
\chi_3 & 2(1+\sqrt{-1})& 1 & 1\cr
\noalign{\hrule}
      }}
\vskip 1pc
Type $\alpha_{1^3}(2)$}}
\vskip 2pc

\hbox{\hskip 1cm
\vbox{
 \vbox{\offinterlineskip\halign{\vrule height9pt depth2.5pt width0pt
         $#$\hfil\quad\vrule\  & $#$\hfil\quad\vrule\ &
         $#$\hfil\quad\vrule\  & $#$\hfil\quad\vrule \cr
          & \hfil f_\chi\hfil & \hfil p=2\hfil & \hfil c\hfil\cr
\noalign{\hrule}
\chi_1 & 4 & 1 & 1\cr
\chi_2 & 2 & 1 & 1\cr
\chi_3 & 4 & 1 & 1\cr
\noalign{\hrule}
      }}
\vskip 1pc
Type $\beta_{1^3}(2)$}
\hskip -10.8cm
\vbox{
 \vbox{\offinterlineskip\halign{\vrule height9pt depth2.5pt width0pt
         $#$\hfil\quad\vrule\  & $#$\hfil\quad\vrule\ &
         $#$\hfil\quad\vrule\  & $#$\hfil\quad\vrule \cr
          & \hfil f_\chi\hfil & \hfil p=3\hfil & \hfil c\hfil\cr
\noalign{\hrule}
\chi_1 & 3 & 1 & 1\cr
\chi_2 & 3 & 1 & 1\cr
\chi_3 & 3 & 1 & 1\cr
\noalign{\hrule}
      }}
\vskip 1pc
Type $\alpha_{1^3}(3)$}
\hskip -10.8cm
\vbox{
 \vbox{\offinterlineskip\halign{\vrule height9pt depth2.5pt width0pt
         $#$\hfil\quad\vrule\  & \hfil$#$\quad\vrule\ &
         $#$\hfil\quad\vrule\  & $#$\hfil\quad\vrule\ \cr
    & \hfil f_\chi\hfil & \hfil p=3\hfil& \hfil c\hfil\cr
\noalign{\hrule}
\chi_1 &   -3\zd& 1 & 1\cr
\chi_2 &       3& 1 & 1\cr
\chi_3 &       3& 1 & 1\cr
\chi_4 & -3\zd^2& 1 & 1\cr
\noalign{\hrule}
      }}
\vskip 1pc
Type $\alpha_{1^4}(3)$}}
\vskip 2pc

\hbox{\hskip 1cm
\vbox{
\vbox{\offinterlineskip\halign{\vrule height9pt depth2.5pt width0pt
         $#$\hfil\ \vrule\  & \hfil$#$\ \vrule\ & $#$\hfil\ \ &
         $#$\hfil\ \vrule & \ $#$\hfil\ & \ $#$\hfil\ \vrule\  \cr
    & \hfil f_\chi\hfil & \multispan2\hfil $p=3$\hfil\vrule& \hfil c\hfil& \cr
\noalign{\hrule}
\chi_1 &   -3\zd& 1 &  & 1& \cr
\chi_2 & -3\zd^2& 1 &  & 1& \cr
\chi_3 &       3&   & 1&  & 1\cr
\chi_4 &       3& 1 & 1& 1& 1\cr
\chi_5 &       3& 1 & 1& 1& 1\cr
\noalign{\hrule}
      }}
\vskip 1pc
Type $\cM(Z_3)$}
\hskip -10.3cm
\vbox{
\vbox{\offinterlineskip\halign{\vrule height9pt depth2.5pt width0pt
         $#$\hfil\ \vrule\  & \hfil$#$\ \vrule\ & $#$\hfil\ \ &
         $#$\hfil\ \vrule & \ $#$\hfil\ & \ $#$\hfil\ \vrule\ \cr
    & \hfil f_\chi\hfil & \multispan2\hfil $p=3$\hfil\vrule& \hfil c\hfil& \cr
\noalign{\hrule}
\chi_1 &   3\zd& 1 &  & 1& \cr
\chi_2 & 3\zd^2& 1 &  & 1& \cr
\chi_3 &      3&   & 1& 1& 1\cr
\chi_4 &      3&   & 1& 1& 1\cr
\chi_5 &      3& 1 & 1& 2& 1\cr
\noalign{\hrule}
      }}
\vskip 1pc
Type $\cM(Z_3)'$}
\hskip -10.5cm
\vbox{
\vbox{\offinterlineskip\halign{\vrule height9pt depth2.5pt width0pt
         $#$\hfil\ \vrule\  & \hfil$#$\ \vrule\ & $#$\hfil\quad\ &
         $#$\hfil\quad\vrule & \ $#$\hfil\ & \ $#$\hfil\  \vrule\cr
    & \hfil f_\chi\hfil & \multispan2\hfil $p=3$\hfil\vrule& \hfil c&\hfil\cr
\noalign{\hrule}
\chi_1 &  3& 1 &  &  1&  \cr
\chi_2 &  3& 1 &  &  1&  \cr
\chi_3 &  3&   & 1&   & 1\cr
\chi_4 &  3&   & 1&   & 1\cr
\chi_5 &  3& 1 & 1&  1& 1\cr
\noalign{\hrule}
      }}
\vskip 1pc
Type $\cM(Z_3)''$}}
\vskip 2pc

\hbox{\hskip 1cm
\vbox{
\vbox{\offinterlineskip\halign{\vrule height9pt depth2.5pt width0pt
         $#$\hfil\quad\vrule\  & $#$\hfil\quad\vrule\  & $#$\hfil\quad &
         $#$\hfil\quad\vrule& \ $#$\hfil\quad& $#$\hfil\quad\vrule \cr
    & \hfil f_\chi\hfil & \multispan2\hfil $p=2$\hfil\vrule& \hfil c\hfil& \cr
\noalign{\hrule}
\chi_1 & 2& 1 &  &1& \cr
\chi_2 & 2&   & 1& &1\cr
\chi_3 & 2& 1 & 1&1&1\cr
\noalign{\hrule}
      }}
\vskip 1pc
Type $\cM(Z_2)$}
\hskip -8.5cm
\vbox{
\vbox{\offinterlineskip\halign{\vrule height10pt depth2.5pt width0pt
         $#$\hfil\quad\vrule\  & \hfil$#$\quad\vrule\ & $#$\hfil\quad &
         $#$\hfil\quad\vrule & \ $#$\hfil\quad & $#$\hfil\quad\vrule &
         \ $#$\hfil\quad\vrule \cr
  & \hfil f_\chi\hfil &\multispan 2 \hfil $p=2$\hfil\vrule &
  \multispan 2 \hfil $p=r$\hfil\vrule & \hfil c \hfil\cr
\noalign{\hrule}
\chi_1 & -2\sqrt{-r}& 1&  &  1&  & 1\cr
\chi_2 &  2\sqrt{-r}&  & 1&  1&  & 1\cr
\chi_3 &           2& 1& 1&   & 1& 2\cr
\noalign{\hrule}  
      }}
\vskip 1pc
Type $\alpha_{G_4}(r)$, $r\ne2$}}
\vskip 2pc

\hbox{\hskip 1cm
\vbox{
\vbox{\offinterlineskip\halign{\vrule height10pt depth2.5pt width0pt
         $#$\hfil\ \vrule\  & $#$\hfil\ \vrule\ & $#$\hfil\quad & 
         $#$\hfil\ \vrule & \ $#$\hfil\quad & $#$\hfil\ \vrule &
         \ $#$\hfil\  \vrule \cr
    & \hfil f_\chi\hfil & \multispan2\hfil $p=2$\hfil\vrule&
   \multispan2\hfil $p=r$\hfil\vrule& \hfil c \hfil \cr
\noalign{\hrule}
 \chi_1& r+\sqrt{4r-15}& 1& & 1& & 1\cr
 \chi_2& r+\sqrt{4r-15}& 1& &  &1& 1\cr
 \chi_3& r-\sqrt{4r-15}&  &1&  &1& 1\cr
 \chi_4& r-\sqrt{4r-15}&  &1& 1& & 1\cr
\noalign{\hrule}
      }}
\vskip 1pc
Type $\alpha_{3\times4}(r)$, $r\ne2$}
\hskip -8cm
\vbox{
\vbox{\offinterlineskip\halign{\vrule height9pt depth2.5pt width0pt
         $#$\hfil\ \vrule\  & $#$\hfil\ \vrule\ & $#$\hfil\ &
         $#$\hfil\ & $#$\hfil\ & $#$\hfil\  \vrule\ & $#$\hfil\ &
         $#$\hfil\ & $#$\hfil\ & $#$\hfil\ \vrule& \ $#$\hfil\ & $#$\hfil\ & $#$\hfil\ \vrule \cr
    & \hfil f_\chi\hfil & \multispan3\hfil $p=2$\hfil& &
      \multispan3\hfil $p=3$\hfil& & & \hfil c \hfil& \cr
\noalign{\hrule}
 \chi_1& 3& 1& & & & 1&  & &  & 1&  & \cr
 \chi_2& 6&  &1& & &  & 1& &  &  & 1& \cr
 \chi_3& 2&  & &1& &  &  &1&  & 1& 1& 1\cr
 \chi_4& 6&  &1&1& & 1&  & & 1& 1&  & 1\cr
 \chi_5& 3&  & & &1&  & 1& & 1&  & 1& 1\cr
\noalign{\hrule}
      }}
\vskip 1pc
Type $\cM(S_3)$}}
\vskip 2pc

\hbox{\hskip 1cm
\vbox{
\vbox{\offinterlineskip\halign{\vrule height9pt depth2.5pt width0pt
         $#$\hfil\ \vrule\  & \hfil$#$\ \vrule\  & $#$\hfil\  &
         $#$\hfil\ & $#$\hfil\ & $#$\hfil\ & $#$\hfil\  \vrule\ &
         $#$\hfil\  & $#$\hfil\ & $#$\hfil\ \vrule& \ $#$\ \hfil\ &
         $#$\hfil\vrule \cr
    & \hfil f_\chi\hfil & \multispan4\hfil $p=2$\hfil& &
      \multispan2\hfil $p=3$\hfil& &\hfil c \hfil& \cr
\noalign{\hrule}
 \chi_1&  6\zd^2& 1& & & & & 1& & & 1& \cr
 \chi_2&    6\zd&  &1& & & & 1& & & 1& \cr
 \chi_3&       3&  & &1& & &  &1& &  & 1\cr
 \chi_4&   -6\zd&  & & &1& & 1&1& & 1& 1\cr
 \chi_5& -6\zd^2&  & & & &1& 1&1& & 1& 1\cr
 \chi_6&  3-\ird& 1& & & &1&  & &1& 2& 1\cr
 \chi_7&  3+\ird&  &1& &1& &  & &1& 2& 1\cr
\noalign{\hrule}
      }}
\vskip 1pc
Type $\alpha_{3\times5}$
}
\hskip -8cm
\vbox{
\vbox{\offinterlineskip\halign{\vrule height9pt depth2.5pt width0pt
         $#$\hfil\ \vrule\  & \hfil$#$\ \vrule\  & $#$\hfil\  &
         $#$\hfil\ & $#$\hfil\ &  $#$\hfil\  \vrule\ \ & $#$\hfil\ \ &
         $#$\hfil\ \ \vrule & \ $#$\hfil\ & $#$\hfil\ \vrule\cr
    & \hfil f_\chi\hfil & \multispan3\hfil $p=2$\hfil& &
     \multispan2\hfil $p=3$\ \hfil\vrule &\hfil c \hfil& \cr
\noalign{\hrule}
 \chi_1&  6\zd^2& 1& & & & 1& &  1& \cr
 \chi_2&    6\zd&  &1& & & 1& &  1& \cr
 \chi_3&   -6\zd&  & &1& & 1& &  1& \cr
 \chi_4& -6\zd^2&  & & &1& 1& &  1& \cr
 \chi_5&  3-\ird& 1& & &1&  &1&  2& \cr
 \chi_6&  3+\ird&  &1&1& &  &1&  2& \cr
\noalign{\hrule}
      }}
\vskip 1pc
Type $\alpha_{3\times5}'$}}
\vskip 2pc

\hbox{\hskip 3cm
\vbox{
\vbox{\offinterlineskip\halign{\vrule height9pt depth2.5pt width0pt
         $#$\hfil\ \vrule\  & \hfil$#$\ \vrule\  & $#$\hfil\quad &
         $#$\hfil\quad& $#$\hfil\quad&
         $#$\hfil\quad& $#$\hfil\quad \vrule\ & $#$\hfil\quad & $#$\hfil\quad&
         $#$\hfil\quad\vrule& \ $#$\hfil\quad& 
         $#$\hfil\ \vrule \cr
    & \hfil f_\chi\hfil & \multispan4\hfil $p=2$\hfil& &
      \multispan2\hfil $p=3$\hfil& & \hfil c \hfil& \cr
\noalign{\hrule}
 \chi_1& -6\zd^2& 1& & & & & 1& & &  1&  \cr
 \chi_2&   -6\zd&  &1& & & & 1& & &  1&  \cr
 \chi_3& -3-\ird&  & &1& & &  &1& &   & 1\cr
 \chi_4& -3+\ird&  & & &1& &  &1& &   & 1\cr
 \chi_5&   -6\zd&  &1&1& & &  & &1&  1& 1\cr
 \chi_6& -6\zd^2& 1& & &1& &  & &1&  1& 1\cr
 \chi_7&       3&  & & & &1& 1& &1&  2& 1\cr
\noalign{\hrule}
      }}
\vskip 1pc
Type $\alpha_{3\times5}''$
}}
\vskip 1pc

Dans la table \ref{typesfamilles}, nous indiquons les types de familles qui
apparaissent dans un groupe de r\'eflexions spetsial exceptionnel donn\'e.

\let\ts=\times
\begin{table}[htbp]\caption{Types de familles} \label{typesfamilles}
\[\begin{array}{|l|ccccccccccccccccc|}
 \cF& 4& 6& 8& 14& 23& 24& 25& 26& 27& 29& 30& 32& 33& 34& 35& 36& 37\cr
\hline
  \alpha_{1^2}(3)& *& *&  & *&  &  & *& *& *&  &  & *& *& *&  &  &  \cr
\alpha_{1^2}(2^2)&  & *& *&  & *& *&  & *&  & *& *&  & *& *&  & *& *\cr
  \alpha_{1^2}(5)&  &  &  &  & *&  &  &  &  &  & *&  &  &  &  &  &  \cr
  \alpha_{1^3}(2)&  &  & *&  &  &  &  &  &  & *&  &  &  &  &  &  &  \cr
   \beta_{1^3}(2)&  &  &  & *&  &  &  &  &  &  &  &  &  &  &  &  &  \cr
  \alpha_{1^3}(3)&  &  &  & *&  &  &  &  & *&  &  & *&  & *&  &  &  \cr
  \alpha_{1^4}(3)&  &  &  &  &  &  &  &  &  &  &  & *&  &  &  &  &  \cr
         \cM(Z_3)&  &  &  &  &  &  &  & *&  &  &  &  &  & *&  &  &  \cr
        \cM(Z_3)'&  &  &  &  &  &  & *&  &  &  &  & *&  &  &  &  &  \cr
       \cM(Z_3)''&  &  &  &  &  &  &  &  &  &  &  &  &  & *&  &  &  \cr
         \cM(Z_2)&  &  &  &  &  &  &  &  & *& *&  &  & *& *& *& *& *\cr
  \alpha_{G_4}(3)& *&  &  &  &  &  &  & *&  &  &  & *&  &  &  &  &  \cr
  \alpha_{G_4}(7)&  &  &  &  &  & *&  &  &  &  &  &  &  &  &  &  &  \cr
\alpha_{3\ts4}(3)&  &  &  &  &  &  &  &  &  &  &  & *&  & *&  &  &  \cr
\alpha_{3\ts4}(5)&  &  &  &  &  &  &  &  & *&  &  &  &  &  &  &  &  \cr
         \cM(S_3)&  &  &  &  &  &  &  &  &  &  &  &  &  &  & *& *& *\cr
   \alpha_{3\ts5}&  &  &  &  &  &  & *&  &  &  &  & *& *& *&  &  &  \cr
  \alpha_{3\ts5}'&  &  &  &  &  &  &  &  &  &  &  & *&  &  &  &  &  \cr
 \alpha_{3\ts5}''&  &  &  &  &  &  &  &  &  &  &  &  &  & *&  &  &  \cr
\hline
\end{array}\]
\end{table}

\subsection{$W=G_4$}

Soit $W=G_4$ le groupe primitif de r\'eflexions complexe de dimension~2 et
d'ordre~24. Pour $\chi\in\Irr(W)$, on pose
$d(\chi)=\chi(1)$ et on note $b(\chi)$ la valuation en~$x$ du degr\'e
fant\^ome $R_\chi$.
Les caract\`eres irr\'eductibles $\chi$ sont
uniquement d\'etermin\'es par les paires $(d(\chi),b(\chi))$ associ\'ees
et nous notons $\chi=\phi_{d,b}$ lorsque $d=d(\chi)$ et $b=b(\chi)$. \par

Le groupe de r\'eflexions $G_4$ a~4 familles de caract\`eres, des types
suivants:

$\bullet$\ d\'efaut $0$~:
$(\phi_{1,0}),(\phi_{3,2})$

$\bullet$ type $\alpha_{1^2}(3)$~:
$(\phi_{2,1},\phi_{2,3})$

$\bullet$ type $\alpha_{G_4}(3)$~:
$(\phi_{1,4},\phi_{1,8},\phi_{2,5})$

\subsection{$W=G_6$}
Les caract\`eres irr\'eductibles de $G_6$ ne sont pas distingu\'es
par leur degr\'e et leur invariant $b$. Dans les cas ambigus, nous
donnons en plus comme troisi\`eme indice le degr\'e en $x$ du degr\'e
fant\^ome $R_\chi$.

Le groupe de r\'eflexions $G_6$ a~5 familles de caract\`eres, des types
suivants:

$\bullet$\ d\'efaut $0$~: $(\phi_{1,0}),(\phi_{3,4})$

$\bullet$ type $\alpha_{1^2}(2^2)$~:
$(\phi_{2,5,13},\phi_{2,7})$

$\bullet$ type $\alpha_{1^2}(3)$~:
$(\phi_{1,10},\phi_{1,14})$

$\bullet$ Autre famille~:

$$ \vbox{\offinterlineskip\halign{\vrule height9.5pt depth3pt width0pt
         $#$\hfil\quad\vrule\ & \hfil$#$\quad\vrule\  & $#$\hfil\quad &
         $#$\hfil\quad & $#$\hfil\quad\vrule\  &
         $#$\hfil\quad & $#$\hfil\quad & $#$\hfil\quad &
         $#$\hfil\quad & $#$\hfil\quad\vrule&
         \ $#$\hfil\quad & $#$\hfil\quad\vrule  \cr
 & \hfil f_\chi\hfil & \multispan2 \hfil $p=2$\hfil& & \multispan4\hfil $p=3$\hfil&
   & c& \cr
\noalign{\hrule}
\phi_{1,4}    &                -4\ird&1& & &1& & &  & & 1& \cr
\phi_{1,8}    &                 4\ird& &1& &1& & &  & & 1& \cr
\phi_{1,6}    &                     4& & &1& &1& &  & &  & 1\cr
\phi_{2,1}    & 2(1+\sqrt{-1})(3+\sqrt{3})&1& &1& & &1&  & & 1& 1\cr
\phi_{2,5,9}  & 2(1+\sqrt{-1})(3-\sqrt{3})& &1&1& & &1&  & & 1& 1\cr
\phi_{2,3,7}  & 2(1-\sqrt{-1})(3-\sqrt{3})&1& &1& & & &1 & & 1& 1\cr
\phi_{2,3,11} & 2(1-\sqrt{-1})(3+\sqrt{3})& &1&1& & & &1 & & 1& 1\cr
\phi_{3,2}    &                     4&1&1&1& & & &  &1& 2& 1\cr
\noalign{\hrule}
      }}
$$

\subsection{$W=G_8$}

Le groupe de r\'eflexions complexe $G_8$ poss\`ede une paire de
caract\`eres irr\'eductibles de m\^emes degr\'e et invariant $b$.
Nous noterons $\phi_{2,7}'$ le conjugu\'e complexe
de $\phi_{2,1}$, et $\phi_{2,7}''$ le conjugu\'e de $\phi_{2,13}$.

Le groupe de r\'eflexions $G_8$ a~5 familles de caract\`eres, des types
suivants:

$\bullet$\ d\'efaut $0$~: $(\phi_{1,0})$

$\bullet$ type $\alpha_{1^2}(2^2)$~:
$(\phi_{4,3},\phi_{4,5})$

$\bullet$ type $\alpha_{1^3}(2)$~:
$(\phi_{2,1},\phi_{2,4},\phi_{2,7}'),(\phi_{3,2},\phi_{3,4},\phi_{3,6})$

$\bullet$ Autre famille~:

$$ \vbox{\offinterlineskip\halign{\vrule height9pt depth3pt width0pt
         $#$\hfil\quad\vrule\ & \hfil$#$\ \vrule\  & $#$\hfil\quad\vrule\  &
         $#$\hfil\quad & $#$\hfil\quad & $#$\hfil\quad &
         $#$\hfil\quad & $#$\hfil\quad & $#$\hfil\quad\vrule &
         \ $#$\hfil\quad\vrule \cr
 & \hfil f_\chi\hfil & \hfil p=2\hfil& & \multispan2 \hfil $p=3$\hfil& & & &c\cr
\noalign{\hrule}
\phi_{1,6}  &         -12&1 &1& & & & &  &1\cr
\phi_{1,18} &         -12&1 & &1& & & &  &1\cr
\phi_{1,12} &           4&1 & & &1& & &  &1\cr
\phi_{2,10} &          12&2 &1&1& & & &  &2\cr
\phi_{2,7}''& -4\sqrt{-1}&2 & & & &1& &  &2\cr
\phi_{2,13} &  4\sqrt{-1}&2 & & & & &1&  &2\cr
\phi_{3,8}  &           4&3 & & & & & &1 &3\cr
\noalign{\hrule}
      }}
$$

Le seul caract\`ere projectif pour $p=2$ dans la famille \`a
7 \'el\'ements
est ind\'ecompo\-sable, car tous les caract\`eres de cette famille se
factorisent (modulo $\fp$) par l'alg\`ebre de groupe $k_\fp[W]$.
Pour $W$, le caract\`ere
projectif correspondant est ind\'ecomposable.

\subsection{$W=G_{14}$}

Les caract\`eres irr\'eductibles de $G_{14}$ ne sont pas distingu\'es
par leur degr\'e et leur invariant $b$. Dans les cas ambigus, nous
donnons en plus comme troisi\`eme indice le degr\'e en $x$ du degr\'e
fant\^ome $R_\chi$.

Le groupe de r\'eflexions $G_{14}$ a~6 familles de caract\`eres, des types
suivants:

$\bullet$\ d\'efaut $0$~: $(\phi_{1,0})$

$\bullet$ type $\beta_{1^3}(2)$~:
$(\phi_{2,9},\phi_{2,12},\phi_{2,15})$

$\bullet$ type $\alpha_{1^2}(3)$~:
$(\phi_{4,5},\phi_{4,7}),(\phi_{1,20},\phi_{1,28})$

$\bullet$ type $\alpha_{1^3}(3)$~:
$(\phi_{3,8},\phi_{3,6,24},\phi_{3,10})$

$\bullet$ Autre famille~:

$$ \vbox{\offinterlineskip\halign{\vrule height9pt depth3pt width0pt
         $#$\hfil\ \vrule\  & \hfil$#$\ \vrule\  & $#$\hfil\quad &
         $#$\hfil\quad & $#$\hfil\quad & $#$\hfil\quad &
         $#$\hfil\quad & $#$\hfil\quad\vrule\  & $#$\hfil\quad &
         $#$\hfil\quad & $#$\hfil\quad & $#$\hfil\quad\vrule&
         \ $#$\hfil\quad & $#$\hfil\quad\vrule \cr
 & \hfil f_\chi\hfil & \multispan5\hfil $p=2$\hfil& & \multispan3\hfil $p=3$\hfil&
  & c&  \cr
\noalign{\hrule}
\phi_{1,8}   &             12\zd^2&1& & & & &  & 1& & & & 1& \cr
\phi_{1,12}  &                   6& &1& & & &  &  &1& & &  & 1\cr
\phi_{1,16}  &               12\zd& & &1& & &  & 1& & & & 1& \cr
\phi_{2,1}   &   12(2-\sqrt{6})\zd& & & &1& &  &  & &1& & 1& 1\cr
\phi_{2,5}   & 12(2-\sqrt{6})\zd^2& & & & &1&  &  & &1& & 1& 1\cr
\phi_{2,7}   &   12(2+\sqrt{6})\zd& & & &1& &  &  & &1& & 1& 1\cr
\phi_{2,11}  & 12(2+\sqrt{6})\zd^2& & & & &1&  &  & &1& & 1& 1\cr
\phi_{2,4}   &              -12\zd& & & &1& &  & 1&1& & & 1& 1\cr
\phi_{2,8}   &            -12\zd^2& & & & &1&  & 1&1& & & 1& 1\cr
\phi_{3,2}   &            -12\zd^2&1& & & &1&  &  & & &1& 2& 1\cr
\phi_{3,4}   &              -12\zd& & &1&1& &  &  & & &1& 2& 1\cr
\phi_{3,6,18}&                   6& & & & & & 1&  & & &1& 2& 1\cr
\phi_{4,3}   &                   6& &1& & & & 1&  & &2& & 2& 2\cr
\noalign{\hrule}
      }}
$$

Avec les m\'ethodes de~\ref{subsec:algo}, il restent deux questions pour la
famille \`a 13
\'el\'ements. Pour $p=2$, on ne trouve pas le deuxi\`eme caract\`ere projectif
mais seulement la somme avec le sixi\`eme. Cependant, un calcul explicite
montre que la repr\'esentation de caract\`ere $\phi_{4,3}$ \`a des facteurs
de degr\'es~3 et~1 modulo~2 et la matrice de d\'ecomposition est comme
indiqu\'ee. Pour $p=3$, il faut montrer que le troisi\`eme projectif $\phi$ est
ind\'ecomposable. On utilise alors le fait que les 4 caract\`eres de
degr\'e~2 intervenant dans $\phi$ se
factorisent par l'alg\`ebre de Hecke $\cH(G_{12})$. Il est alors
facile de voir qu'ils ont m\^eme r\'eduction modulo~3.

\subsection{$W=G_{23}=W(H_3)$}

Le groupe de r\'eflexions $G_{23}$ a~7 familles de caract\`eres, des types
suivants:

$\bullet$\ d\'efaut $0$~: $(\phi_{1,0}),(\phi_{5,2}),(\phi_{5,5}),
(\phi_{1,15})$

$\bullet$ type $\alpha_{1^2}(2^2)$~: $(\phi_{4,3},\phi_{4,4})$

$\bullet$ type $\alpha_{1^2}(5)$~: $(\phi_{3,1},\phi_{3,3})$,
$(\phi_{3,6},\phi_{3,8})$

\subsection{$W=G_{24}$}

Le groupe de r\'eflexions $G_{24}$ a~7 familles de caract\`eres, des types
suivants:

$\bullet$\ d\'efaut $0$~: $(\phi_{1,0}),(\phi_{7,3}),
(\phi_{7,6}),(\phi_{1,21})$

$\bullet$ type $\alpha_{1^2}(2^2)$~:
$(\phi_{8,4},\phi_{8,5})$

$\bullet$ type $\alpha_{G_4}(7)$~:
$(\phi_{3,3},\phi_{3,1},\phi_{6,2})$
$(\phi_{3,8},\phi_{3,10},\phi_{6,9})$

Pour \'etablir ce resultat il faut voir que le projectif
$\phi_{3,1}+\phi_{3,3}+2\phi_{6,2}$ n'est pas ind\'ecomposable modulo~2
(et de m\^eme pour
le produit avec le caract\`ere d\'eterminant).
Puisque $\phi_{3,1}$ et $\phi_{3,3}$ ont
des r\'eductions diff\'erentes pour l'alg\`ebre $k_\fp[W]$, alors il en
est de m\^eme pour $\cH$.

\subsection{$W=G_{25}$}

Pour distinguer les caract\`eres de $W=G_{25}$
ayant m\^emes valeurs $d,b$ nous
choissisons la notation telle que $(\phi_{3,1},\phi_{3,5}')$,
$(\phi_{3,5}'',\phi_{3,13}'')$, $(\phi_{3,13}',\phi_{3,17})$,
$(\phi_{6,2},\phi_{6,4}')$, $(\phi_{6,4}'',\phi_{6,8}'')$ et
$(\phi_{6,8}',\phi_{6,10})$ sont des paires de caract\`eres
conjugu\'es.

Le groupe de r\'eflexions $G_{25}$ a~7 familles de caract\`eres, des types
suivants:

$\bullet$\ d\'efaut $0$~: $(\phi_{1,0})$

$\bullet$ type $\alpha_{1^2}(3)$~:
$(\phi_{3,1},\phi_{3,5}'),(\phi_{8,6},\phi_{8,9}),(\phi_{6,8}',\phi_{6,10})$

$\bullet$ type $\cM(Z_3)'$~: $(\phi_{2,3},\phi_{2,9},\phi_{6,2},\phi_{6,4}',
  \phi_{8,3})$

$\bullet$ type $\alpha_{3\times5}$~:
$(\phi_{3,5}'',\phi_{3,13}'',\phi_{3,6},\phi_{6,4}'',\phi_{6,8}'',\phi_{9,5},
  \phi_{9,7})$

$\bullet$ Autre famille~:

$$ \vbox{\offinterlineskip\halign{\vrule height9pt depth3pt width0pt
         $#$\hfil\ \vrule\  & \hfil$#$\ \vrule\ & $#$\hfil\quad &
         $#$\hfil\quad & $#$\hfil\quad & $#$\hfil\quad\vrule\  &
         $#$\hfil\quad & $#$\hfil\quad\vrule& \ $#$\hfil\quad\vrule \cr
 & \hfil f_\chi\hfil & \multispan3 \hfil $p=2$\hfil& &
    \multispan{2} \hfil $p=3$\hfil\vrule& c \cr
\noalign{\hrule}
\phi_{1,12} &      -6&1& & & &  1& & 1\cr
\phi_{1,24} &      -6& &1& & &  1& & 1\cr
\phi_{2,15} &       6&1&1& & &   &1& 2\cr
\phi_{3,13}'&   -3\zd& & &1& &  1&1& 3\cr
\phi_{3,17} & -3\zd^2& & & &1&  1&1& 3\cr
\noalign{\hrule}
      }}
$$

Plusieurs points sont \`a justifier. Tout d'abord, il faut d\'emontrer
qu'il y a un module simple de dimension $6$ dans la famille
$\cM(Z_3)'$ pour $p=3$.

Supposons $p=3$. Soit $C$ le bloc de $\cO_\fp\cH$ correspondant \`a la famille
de type $\cM(Z_3)'$. Soit $W'$ un sous-groupe parabolique de $W$
de type $G(3,1,1)\times G(3,1,1)$, de g\'en\'erateurs $s,u$ (conjugu\'es dans
$W$) et $\cH'$ son alg\`ebre de Hecke.

Le lemme \ref{projrela} montre que
l'alg\`ebre $C$ est projective relativement \`a $\cO_\fp\cH'$.

Soit $S$ un $\cO_\fp\cH$-module libre sur $\cO_\fp$ tel que
$K\otimes_{\cO_\fp} S$ a pour caract\`ere $\phi_{2,3}$.
Soit $\bar{S}$ sa r\'eduction modulo $3$. C'est un
$k_\fp\cH$-module simple.

Soit $T$ le $C$-module simple non isomorphe \`a $\bar{S}$. Nous allons
montrer que $T$ est de dimension $6$.

On note $\psi_{i,j}$ le caract\`ere $s\mapsto \det(s)^i,\ u\mapsto
\det(t)^j$ de $W'$.

Soit $Y=C\otimes_{\cO_\fp\cH'} L$, o\`u $L$ est le $\cO_\fp\cH'$-module
libre de rang $1$ sur $\cO_\fp$ tel que $K\otimes_{\cO_\fp} L$ est
de caract\`ere $\psi_{1,1}$.
Alors, le caract\`ere de $Y$ est $\phi_{2,3}+\phi_{6,2}$.
On a $\Res_{\cH'}S\iso L\oplus P$, o\`u $P$ est le
$\cO_\fp\cH'$-module projectif de caract\`ere $\phi_{0,0}$.

Puisque $C$ est projective relativement \`a $\cO_\fp\cH'$,  la
surjection canonique
$$C\otimes_{\cO_\fp\cH'}\Res_{\cH'}S\twoheadrightarrow S$$
est scind\'ee. On en d\'eduit que $S|Y\oplus\text{ projectif}$, donc
que $S|Y$. Soit $Y'$ tel que $Y\simeq S\oplus Y'$.
On a $\Hom(\bar{Y},\bar{S})\simeq \Hom(\bar{L},\Res_{\cH'}\bar{S})\simeq k_\fp$.
On en d\'eduit que $\Hom(\bar{Y}',\bar{S})=0$. De m\^eme,
$\Hom(\bar{S},\bar{Y}')=0$. Puisque la multiplicit\'e de $T$ dans
$\bar{Y}'$ est $1$, on en d\'eduit que $\bar{Y}'=T$. Ceci
montre que $T$ est de dimension $6$.

\smallskip
Pour achever la d\'emonstration des r\'esultats annonc\'es plus haut, il
reste \`a prouver que
$\phi_{9,5}$ et $\phi_{9,7}$ ne sont pas dans le m\^eme $3$-bloc que
$\phi_{3,5}'',\phi_{3,13}'',\phi_{3,6},\phi_{6,4}''$ et $\phi_{6,8}''$.
Nous le d\'eduisons (cf. lemme \ref{projrela}) du fait que
$\Tr_{\cH'}^{\cH}(1)$ agit par un scalaire inversible de
$\cO_\fp$ sur les $5$ derniers caract\`eres, mais pas sur
$\phi_{9,5}$ et $\phi_{9,7}$.

\subsection{$W=G_{26}$}

Pour distinguer les caract\`eres de $W=G_{26}$ ayant m\^emes $d,b$,
nous choissisons la notation telle que $(\phi_{3,1},\phi_{3,5}')$,
$(\phi_{3,5}'',\phi_{3,13}'')$, $(\phi_{3,13}',\phi_{3,17})$,
$(\phi_{3,4},\phi_{3,8}')$, $(\phi_{3,8}'',\phi_{3,16}''),\!$
$(\phi_{3,16}',\phi_{3,20})$,
$(\phi_{6,2},\phi_{6,4}'),(\phi_{6,4}'',\phi_{6,8}'')$,
$(\phi_{6,8}',\phi_{6,10}),(\phi_{6,5},\phi_{6,7}'),
(\phi_{6,7}'',\phi_{6,11}'')$ et $(\phi_{6,11}',\phi_{6,13})$
sont des paires de caract\`eres conjugu\'es. De plus,
$R_\chi\equiv x^{24}\pmod{x^{23}}$ pour $\chi=\phi_{8,6}'$ et
$R_\chi\equiv x^{27}\pmod{x^{26}}$ pour $\chi=\phi_{8,9}'$.

Le groupe de r\'eflexions $G_{26}$ a~10 familles de caract\`eres, des types
suivants:

$\bullet$\ d\'efaut $0$~: $(\phi_{1,0})$

$\bullet$ type $\alpha_{1^2}(2^2)$~:
$(\phi_{8,3}, \phi_{8,6}')$

$\bullet$ type $\alpha_{1^2}(3)$~:
$(\phi_{6,5},\phi_{6,7}')$, $(\phi_{3,16}',\phi_{3,20})$

$\bullet$ type $\alpha_{G_4}(3)$~:
$(\phi_{1,21},\phi_{1,33},\phi_{2,24})$

$\bullet$ type $\cM(Z_3)$~:
$(\phi_{2,3},\phi_{2,9},\phi_{1,9},\phi_{3,1},\phi_{3,5}')$,
$(\phi_{3,4},\phi_{3,8}',\phi_{3,6},\phi_{6,4}',\phi_{6,2})$,
$\!(\phi_{3,13}',\phi_{3,17},\phi_{3,15},\!$ $\phi_{6,13},\phi_{6,11}')$

$\bullet$ Autres familles~:

$$ \vbox{\offinterlineskip\halign{\vrule height8.5pt depth2.5pt width0pt
         $#$\hfil\ \vrule\  & \hfil$#$\ \vrule\  & $#$\hfil\  &
         $#$\hfil\  & $#$\hfil\  & $#$\hfil\ \vrule\  &
         $#$\hfil\  & $#$\hfil\  & $#$\hfil\ \vrule  &
         \ $#$\hfil\  & $#$\hfil\ \vrule \cr
 & \hfil f_\chi\hfil & \multispan{3} \hfil $p=2$\hfil& &
   \multispan{2} \hfil $p=3$\hfil& &c& \cr
\noalign{\hrule}
 \phi_{3,5}''& -3-\ird&1& & &  &1& & & 1& \cr
\phi_{3,13}''& -3+\ird& &1& &  &1& & & 1& \cr
 \phi_{6,4}''&  3-\ird& & &1&  & &1& & 1& 1\cr
 \phi_{6,8}''&  3+\ird& & & &1 & &1& & 1& 1\cr
 \phi_{9,7}  &  3-\ird& &1&1&  & & &1& 2& 1\cr
 \phi_{9,5}  &  3+\ird&1& & &1 & & &1& 2& 1\cr
\noalign{\hrule}
      }}
$$

$$ \vbox{\offinterlineskip\halign{\vrule height9pt depth2.5pt width0pt
         $#$\hfil\ \vrule\  & \hfil$#$\ \vrule\ &
         $#$\hfil\  & $#$\hfil\  & $#$\hfil\  &
         $#$\hfil\  & $#$\hfil\  & $#$\hfil\  & $#$\hfil\  &
         $#$\hfil\ \vrule\  & $#$\hfil\  & $#$\hfil\  &
         $#$\hfil\  & $#$\hfil\  & $#$\hfil\ \vrule &
         \ $#$\hfil\  & $#$\hfil\  & $#$\hfil\ \vrule \cr
 & \hfil f_\chi\hfil& \multispan{7} \hfil $p=2$\hfil& &
   \multispan{4} \hfil $p=3$\hfil& & & c& \cr
\noalign{\hrule}
 \phi_{1,12} &      18\ird&1& & & & & & &  &1& & & & & 1& & \cr
 \phi_{1,24} &     -18\ird& &1& & & & & &  &1& & & & & 1& & \cr
 \phi_{2,15} &           6&1&1& & & & & &  & &1& & & & 2& & \cr
 \phi_{2,18} &       9\ird& & & & & & & &  & & &1& & &  & 1& \cr
 \phi_{2,12} &      -9\ird& & & & & & & &  & & &1& & &  & 1& \cr
 \phi_{3,8}''&   3(3-\ird)& & &1& & & & &  &1& &1& & & 1& 1& \cr
\phi_{3,16}''&   3(3+\ird)& & & &1& & & &  &1& &1& & & 1& 1& \cr
 \phi_{6,8}' &    3\ird\zd& & & & & & & &  &1& & &1& & 3& & 1\cr
 \phi_{6,10} & -3\ird\zd^2& & & & & & & &  &1& & &1& & 3& & 1\cr
\phi_{6,11}''&     -6\zd^2& & & & &1& & &  & & & & &1& 1& 1& 1\cr
 \phi_{6,7}''&       -6\zd& & & & & &1& &  & & & & &1& 1& 1& 1\cr
 \phi_{8,9}' &           6& & & & & & &1&  & &1& & &1& 3& 1& 1\cr
 \phi_{8,9}''&           6& & & & & & & &1 & &1& & &1& 3& 1& 1\cr
 \phi_{8,12} &      18\ird& & & & & & & &1 &1& &1&1& & 3& 1& 1\cr
 \phi_{8,6}''&     -18\ird& & & & & & &1&  &1& &1&1& & 3& 1& 1\cr
 \phi_{9,8}  &   3(3-\ird)& & &1& &1& & &  & & &2&1& & 2& 2& 1\cr
 \phi_{9,10} &   3(3+\ird)& & & &1& &1& &  & & &2&1& & 2& 2& 1\cr
\noalign{\hrule}
      }}
$$
Dans la table nous avons omis les projectifs appartenant \`a des blocs de
$p$-d\'efaut~0.

\subsection{$W=G_{27}$}

Le groupe de r\'eflexions complexe $G_{27}$ a cinq paires de caract\`eres
de m\^emes degr\'e et invariant $b$. Nous distinguons ceux de degr\'e~3 et~8
en donnant comme troisi\`eme indice le degr\'e en $x$ du degr\'e fant\^ome.
\par
Le produit par le caract\`ere d\'eterminant $\epsilon$ donne un
automorphisme des matrices de d\'ecomposition. Nous utilisons ce fait en
donnant le r\'esultat pour $\chi$ et $\chi\otimes\epsilon$ dans la m\^eme ligne.

Le groupe de r\'eflexions $G_{27}$ a~11 familles de caract\`eres, des types
suivants:

$\bullet$\ d\'efaut $0$~: $(\phi_{1,0}),(\phi_{1,45})$

$\bullet$ type $\alpha_{1^2}(3)$~:
$(\phi_{15,7},\phi_{15,5})$, $(\phi_{15,8},\phi_{15,10})$

$\bullet$ type $\cM(Z_2)$~:
$(\phi_{5,6}',\phi_{5,6}'',\phi_{10,3})$,
$(\phi_{5,15}'',\phi_{5,15}',\phi_{10,12})$

$\bullet$ type $\alpha_{1^3}(3)$~:
$(\phi_{9,4},\phi_{9,8},\phi_{9,6})$, $(\phi_{9,9},\phi_{9,11},\phi_{9,13})$

$\bullet$ type $\alpha_{3\times4}(5)$~:
$(\phi_{8,6},\phi_{8,9,39},\phi_{8,9,33},\phi_{8,12})$

$\bullet$ Autre famille~:

$$ \vbox{\offinterlineskip\halign{\vrule height9pt depth3pt width0pt
         $#$\hfil\quad\vrule\  & $#$\hfil\quad\vrule\  & \hfil$#$\quad\vrule\ 
         & $#$\hfil\quad &
         $#$\hfil\quad & $#$\hfil\quad & $#$\hfil\quad\vrule\  & 
         $#$\hfil\quad & $#$\hfil\quad & $#$\hfil\quad\vrule\  &
         $#$\hfil\quad & $#$\hfil\quad\vrule & \ $#$\hfil\quad\vrule \cr
 \chi& \chi\otimes\eps& f_\chi& \multispan3 \hfil $p=2$\hfil& &
   \multispan{2} \hfil $p=3$\hfil& & \multispan{1} \hfil $p=5$\hfil& & c\cr
\noalign{\hrule}
\phi_{3,1}   & \phi_{3,20,44}& -2\sqrt{-15}\zd^2&1& & & & 1& & & 1& & 1\cr
\phi_{3,5,29}& \phi_{3,16}   &    2\sqrt{-15}\zd& &1& & & 1& & &  &1& 1\cr
\phi_{3,5,23}& \phi_{3,22}   &   -2\sqrt{-15}\zd& & &1& &  &1& &  &1& 1\cr
\phi_{3,7}   & \phi_{3,20,38}&  2\sqrt{-15}\zd^2& & & &1&  &1& & 1& & 1\cr
\phi_{6,2}   & \phi_{6,19}   &         -2\ird\zd& &1&1& &  & &1&  & & 2\cr
\phi_{6,4}   & \phi_{6,17}   &        2\ird\zd^2&1& & &1&  & &1&  & & 2\cr
\noalign{\hrule}
      }}
$$

\subsection{$W=G_{28}=W(F_4)$}

Les r\'esultats ont \'et\'e obtenus par Gyoja \cite{Gy}. Pour les
caract\`eres de $W(F_4)$, nous utilisons la notation standard
(cf \cite[Table~C.3]{GP}). 

Il y a~11 familles de caract\`eres, des types suivants:

$\bullet$ type $\cM(Z_2)$~: $(\phi_{4,1},\phi_{2,4}'',\phi_{2,4}')$,
  $(\phi_{4,13},\phi_{2,16}',\phi_{2,16}'')$

$\bullet$ Autre famille~:

$$ \vbox{\offinterlineskip\halign{\vrule height9pt depth3pt width0pt
         $#$\hfil\quad\vrule\  & \hfil$#$\quad\vrule\ &  $#$\hfil\quad &
         $#$\hfil\quad & $#$\hfil\quad & $#$\hfil\quad &
         $#$\hfil\quad\vrule\  & $#$\hfil\quad & $#$\hfil\quad\vrule\  &
         $#$\hfil\quad & $#$\hfil\quad & $#$\hfil\quad & $#$\hfil\quad &
         $#$\hfil\quad\vrule\cr
 \chi& f_\chi& \multispan4 \hfil $p=2$\hfil& & \multispan{1}
   \hfil $p=3$\hfil& & & & c& & \cr
\noalign{\hrule}
\phi_{1,12}''& 8& 1& & & & &   & &   1& & & & \cr
\phi_{1,12}' & 8&  &1& & & &   & &    &1& & & \cr
\phi_{4,7}'' & 4&  & &1& & &   & &   1& &1& & \cr
\phi_{4,7}'  & 4&  & & &1& &   & &    &1& &1& \cr
\phi_{4,8}   & 8&  & & & &1&   & &    & & & &1\cr
\phi_{6,6}'  & 3&  & & & & &  1& &    & &1&1& \cr
\phi_{6,6}'' &12& 1&1& & &1&   &1&   1&1& & &1\cr
\phi_{9,6}'' & 8& 1& &1& &1&   & &   2& &1& &1\cr
\phi_{9,6}'  & 8&  &1& &1&1&   & &    &2& &1&1\cr
\phi_{12,4}  &24&  & &1&1&1&  1&1&   1&1&1&1&1\cr
\phi_{16,5}  & 4&  & &1&1&2&   & &   1&1&1&1&2\cr
\noalign{\hrule}
      }}
$$

Tous les autres familles sont de d\'efaut $0$.

\subsection{$W=G_{29}$}

Pour $W=G_{29}$, nous notons $\phi_{6,10}'$, $\phi_{6,10}''$ les deux
caract\`eres r\'eels de degr\'e~6 et d'invariant~$b$ \'egal \`a~$10$. En outre,
$\phi_{15,4}''$ est le caract\`ere de degr\'e~15 intervenant dans
$\phi_{4,3}\otimes\phi_{4,1}$ et $\phi_{15,12}''=\phi_{15,4}''\otimes\eps$
avec $\eps=\phi_{1,40}$ le caract\`ere d\'eterminant.

Le groupe de r\'eflexions $G_{29}$ a~15 familles de caract\`eres, des types
suivants:

$\bullet$\ d\'efaut $0$~: $(\phi_{1,0}),(\phi_{10,2}),(\phi_{15,4}'),
(\phi_{15,12}'),(\phi_{10,18}),(\phi_{1,40})$

$\bullet$ type $\alpha_{1^3}(2)$:
$\!(\phi_{4,3},\phi_{4,4},\phi_{4,1})$,
$\!\!(\phi_{20,5},\phi_{20,6},\phi_{20,7})$,
$\!\!(\phi_{20,11},\phi_{20,10},\phi_{20,9})$,
$\!\!(\phi_{4,21},$ $\phi_{4,24},\phi_{4,23})\!\!$

$\bullet$ type $\alpha_{1^2}(2^2)$~:
$(\phi_{16,3},\phi_{16,5})$, $(\phi_{16,13},\phi_{16,15})$

$\bullet$ type $\cM(Z_2)$~:
$(\phi_{5,8},\phi_{10,6},\phi_{15,4}'')$,
$(\phi_{5,16},\phi_{10,14},\phi_{15,12}'')$

$\bullet$ Autre famille~:

$$ \vbox{\offinterlineskip\halign{\vrule height9pt depth3pt width0pt
         $#$\hfil\quad\vrule\  & \hfil$#$\quad\vrule\ &  $#$\hfil\quad &
         $#$\hfil\quad & $#$\hfil\quad &
         $#$\hfil\quad\vrule\  & $#$\hfil\quad & $#$\hfil\quad &
         $#$\hfil\quad & $#$\hfil\ \vrule\  &
         $#$\hfil\quad & $#$\hfil\quad & $#$\hfil\ \vrule \cr
 \chi& f_\chi& \multispan3 \hfil $p=2$\hfil& & \multispan{3}
   \hfil $p=5$\hfil& & & c& \cr
\noalign{\hrule}
\phi_{6,10}'  &   5& 1& & & &  1& & & &   1& & \cr
\phi_{6,10}''  &  5&  &1& & &   &1& & &    &1& \cr
\phi_{6,10}''' &-20&  & &1& &   & &1& &    & &1\cr
\phi_{6,10}''''&-20&  & &1& &   & & &1&    & &1\cr
\phi_{6,12}  &    4&  & &1& &   & & & &    & &1\cr
\phi_{24,6}  &   20&  & &2&1&  1&1&1&1&   1&1&2\cr
\phi_{24,7}  &    4&  & &2&1&   & & & &   1&1&2\cr
\phi_{24,9}  &    4&  & &2&1&   & & & &   1&1&2\cr
\phi_{30,8}  &    4&  & &3&1&   & & & &   1&1&3\cr
\noalign{\hrule}
   ??        &     & & &*& &   & & &  & & & \cr
\noalign{\hrule}
      }}
$$

Pour $W=G_{29}$ nous ne savons pas d\'emontrer que le troisi\`eme caract\`ere
projectif modulo~2 dans la famille \`a~9 \'el\'ements est ind\'ecomposable.

Nous l'indiquons par l'$*$ sous la colonne correspondante de la table.

\subsection{$W=G_{30}=W(H_4)$}

Pour $W=G_{30}$ le caract\`ere $\phi_{30,10}'$ de degr\'e~30 est celui
pour lequel
$\phi_{30,10}'(c)=\phi_{4,1}(c)$ pour un \'el\'ement de Coxeter $c\in W$,
tandis que $\phi_{30,10}''(c)=\phi_{4,7}(c)$.

Le groupe de r\'eflexions $G_{30}$ a~13 familles de caract\`eres, des types
suivants:

$\bullet$\ d\'efaut $0$~: $(\phi_{1,0}),(\phi_{25,4}),(\phi_{36,5}),
(\phi_{36,15}),(\phi_{25,16}),(\phi_{1,60})$

$\bullet$ type $\alpha_{1^2}(5)$~: $(\phi_{4,1},\phi_{4,7})$,
$(\phi_{9,2},\phi_{9,6})$, $(\phi_{9,22},\phi_{9,26})$,
$(\phi_{4,31},\phi_{4,37})$

$\bullet$ type $\alpha_{1^2}(2^2)$~: $(\phi_{16,3},\phi_{16,6})$,
$(\phi_{16,18},\phi_{16,21})$

$\bullet$ Autre famille~:

$$ \vbox{\offinterlineskip\halign{\vrule height9pt depth3pt width0pt
         $#$\hfil\ \vrule\  & \hfil$#$\ \vrule\  & $#$\hfil\  &
         $#$\hfil\  & $#$\hfil\  & $#$\hfil\  & $#$\hfil\  &
         $#$\hfil\  & $#$\hfil\  & $#$\hfil\  &$#$\hfil\ \vrule\ &
         $#$\hfil\  & $#$\hfil\ & $#$\hfil\ & $#$\hfil\ &
         $#$\hfil\ & $#$\hfil\ \vrule\  &
         $#$\hfil\ & $#$\hfil\ & $#$\hfil\ & $#$\hfil\ &
         $#$\hfil\ & $#$\hfil\ \vrule\  &
         \ $#$\hfil\ & \ $#$\hfil\ & \ $#$\hfil\ \vrule \cr
 \chi& f_\chi& \multispan8 \hfil $p=2$\hfil& & \multispan{5}
   \hfil $p=3$\hfil& & \multispan{5}  \hfil $p=5$\hfil& & & c&  \cr
\noalign{\hrule}
 \phi_{6,12}   &15(3+\sqrt5)&1& & & & & & & & &  1& & & & & &  1& & & & & &   1&  & \cr
 \phi_{6,20}   &15(3-\sqrt5)& &1& & & & & & & &   &1& & & & &  1& & & & & &   1&  & \cr
 \phi_{8,13}   &          10& & &1& & & & & & &   & & & & & &   &1& & & & &    & 1& \cr
 \phi_{10,12}  &          10& & & &1& & & & & &   & & & & & &   & &1& & & &    &  & 1\cr
 \phi_{18,10}  &          10& & &1&1& & & & & &   & & & & & &   & & &1& & &    & 1& 1\cr
 \phi_{24,7}   &15(7+3\sqrt5)& & & & &1& & & & &   & &1& & & &   &1& & &1& &  1& 1& 1\cr
 \phi_{24,11}  &15(7-3\sqrt5)& & & & & &1& & & &   & & &1& & &   &1& & &1& &  1& 1& 1\cr
 \phi_{30,10}' &15(3-\sqrt5)& &1& & &1& & & & &  1& & & &1& &   & & & & &1&   2& 1& 1\cr
 \phi_{30,10}''&15(3+\sqrt5)&1& & & & &1& & & &   &1& & & &1&   & & & & &1&   2& 1& 1\cr
 \phi_{8,12}   &           8& & & & & & &1& & &   & & & & & &   & & & & & &    & 1& \cr
 \phi_{16,11}  &          20& & & & & & & &1& &   & & & & & &   & & & &1& &   1&  & 1\cr
 \phi_{16,13}  &          20& & & & & & & & &1&   & & & & & &   & & & &1& &   1&  & 1\cr
 \phi_{24,6}   &120(9+4\sqrt5)& & & & & & &1&1& &   & & & & &1& 1& & &1& & &  1& 1& 1\cr
 \phi_{24,12}  &120(9-4\sqrt5)& & & & & & &1& &1&   & & & &1& & 1& & &1& & &  1& 1& 1\cr
 \phi_{40,8}   &          40& & & & & & &1&1&1&   & & & & & &   & &1& & &1&   2& 1& 2\cr
 \phi_{48,9}   &          12& & & & & & &2&1&1&   & &1&1& & &   & & & & & &   2& 2& 2\cr
\noalign{\hrule}
      }}
$$

En utilisant le logiciel {\sf MeatAxe}, J\"urgen M\"uller a d\'emontr\'e
que les caract\`eres de d\'egre~30 restent
irr\'eductibles modulo~5. \par
Les caract\`eres constructibles de $W(H_4)$ sont d\'ej\`a dans
Alvis--Lusztig \cite{AlLu}.

\subsection{$W=G_{32}$}

Pour distinguer les caract\`eres de $W=G_{32}$ ayant m\^eme valeurs $d,b$,
nous choisissons la notation telle que $(\phi_{20,3},\phi_{20,9}')$,
$(\phi_{20,9}'',\phi_{20,21})$, $(\phi_{20,13},\phi_{20,29}')$,
$(\phi_{20,29}'',\phi_{20,31})$, $(\phi_{30,12}',\phi_{30,24})$,
$(\phi_{30,16},\phi_{30,20}')$, $(\phi_{30,20}'',\phi_{30,28})$,
$(\phi_{60,9},\phi_{60,15}')$, $(\phi_{60,7},\phi_{60,11}')$ et
$(\phi_{60,11}'',\phi_{60,13})$ sont des paires de caract\`eres
conjugu\'es.

Le groupe de r\'eflexions $G_{32}$ a~16 familles de caract\`eres, des types
suivants:

$\bullet$\ d\'efaut $0$~: $(\phi_{1,0})$

$\bullet$ type $\alpha_{1^2}(3)$~: $(\phi_{4,1},\phi_{4,11})$,
$(\phi_{20,3},\phi_{20,9}')$,

$\bullet$ type $\alpha_{G_4}(3)$~: $(\phi_{10,30},\phi_{10,42},\phi_{20,33})$

$\bullet$ type $\alpha_{3\times5}$~:
$(\phi_{6,8},\phi_{6,28},\phi_{24,6},\phi_{30,4},\phi_{30,8},\phi_{36,5},
  \phi_{36,7})$,
$(\phi_{15,22},\phi_{15,38},\phi_{15,24}, \phi_{30,20}''$,\break
  $\phi_{30,28},\phi_{45,22},\phi_{45,26})$

$\bullet$ type $\alpha_{3\times5}'$~:
$(\phi_{20,5},\phi_{20,19},\phi_{20,7},\phi_{20,17},\phi_{40,8},\phi_{40,10})$,

$\bullet$ type $\alpha_{3\times4}(3)$~:
$(\phi_{64,13},\phi_{64,16},\phi_{64,8},\phi_{64,11})$

$\bullet$ type $\cM(Z_3)'$~:
$(\phi_{5,20},\phi_{5,4},\phi_{10,2},\phi_{10,10},\phi_{15,6})$

$\bullet$ type $\alpha_{1^4}(3)$~:
 $(\phi_{20,25},\phi_{20,29}'',\phi_{20,31},\phi_{20,35})$,
 $(\phi_{60,9},\phi_{60,11}'',\phi_{60,13},\phi_{60,15}')$

$\bullet$ type $\alpha_{1^3}(3)$~: $(\phi_{81,10},\phi_{81,12},\phi_{81,14})$

$\bullet$ Autres familles~:

$$\vbox{\offinterlineskip\halign{\vrule height9pt depth2pt width0pt
         $#$\hfil\ \vrule\  & \hfil$#$\ \vrule\  & \ $#$\hfil\  &
         \ $#$\hfil\  & \ $#$\hfil\  & \ $#$\hfil\  & \ $#$\hfil\  &
         \ $#$\hfil\  & \ $#$\hfil\  & \ $#$\hfil\  & \ $#$\hfil\ \vrule\ 
         & \ $#$\hfil\  & \ $#$\hfil\  & \ $#$\hfil\  & \ $#$\hfil\  &
         \ $#$\hfil\ \vrule\   & \ $#$\hfil\  &
         \ $#$\hfil\  & \ $#$\hfil\ \vrule\cr
 \chi& f_\chi& \multispan{8} \hfil $p=2$\hfil& & \multispan{4}
   \hfil $p=3$\hfil& & & c& \cr
\noalign{\hrule}
 \phi_{5,12}  &      24&1& & & & & & & & &  1& & & & &  1&  &  \cr
 \phi_{5,36}  &      24& &1& & & & & & & &  1& & & & &  1&  &  \cr
 \phi_{15,8}  &  6\zd^2& & &1& & & & & & &   &1& & & &  1& 1&  \cr
 \phi_{15,16} &    6\zd& & & &1& & & & & &   &1& & & &  1& 1&  \cr
 \phi_{20,12} &      24& & & & &1& & & & &  1&1& & & &  2& 1&  \cr
 \phi_{20,9}''&      12& & & & & &1& & & &   & &1& & &  1&  & 1\cr
 \phi_{20,21} &      12& & & & & & &1& & &   & &1& & &  1&  & 1\cr
\phi_{30,12}''&       4&1&1& & &1& & & & &   & & & & &  4& 1&  \cr
 \phi_{45,6}  &      24&1& & & &1&1& & & &   & & &1& &  4& 1& 1\cr
 \phi_{45,18} &      24& &1& & &1& &1& & &   & & &1& &  4& 1& 1\cr
 \phi_{45,14} & -6\zd^2& & & & & & & &1& &   & & &1& &  4& 1& 1\cr
 \phi_{45,10} &   -6\zd& & & & & & & & &1&   & & &1& &  4& 1& 1\cr
 \phi_{60,12} &       8& & & & &1&1&1& & &   & & & & &  4& 1& 2\cr
 \phi_{60,11}'& -6\zd^2& & &1& & & & &1& &   & &1& &1&  5& 2& 1\cr
 \phi_{60,7}  &   -6\zd& & & &1& & & & &1&   & &1& &1&  5& 2& 1\cr
 \phi_{80,9}  &      12& & & & &2&1&1& & &   & &2& &1&  6& 2& 2\cr
\noalign{\hrule}
 \    ??      &        & & & & & & & & & &   & &*& & &  & & \cr
\noalign{\hrule}
      }}$$

$$ \vbox{\offinterlineskip\halign{\vrule height9pt depth3pt width0pt
         $#$\hfil\ \vrule\  & \hfil$#$\ \vrule\  & $#$\hfil\  &
         $#$\hfil\  & $#$\hfil\  & $#$\hfil\  & $#$\hfil\  &
         $#$\hfil\  & $#$\hfil\  & $#$\hfil\  & $#$\hfil\  &
         $#$\hfil\ \vrule\  & $#$\hfil\  & $#$\hfil\  &
         $#$\hfil\  & $#$\hfil\  & $#$\hfil\ \vrule\  &
         \ $#$\hfil\  & \ $#$\hfil\  & \ $#$\hfil\ \vrule\cr
 \chi& f_\chi& \multispan{9} \hfil $p=2$\hfil& & \multispan{4}
   \hfil $p=3$\hfil& & &c& \cr
\noalign{\hrule}
 \phi_{10,34} & -6\ird&1& & & & & & & & & & 1& & & & &  1&  &  \cr
 \phi_{20,13} &     -6& &1& & & & & & & & &  &1& & & &  1& 1&  \cr
 \phi_{20,16} &  6\ird& & &1& & & & & & & &  & &1& & &   & 1& 1\cr
 \phi_{30,16} & -6\ird&1&1& & & & & & & & &  & & &1& &  2& 1&  \cr
 \phi_{60,16} & -6\ird& & & &1& & & & & & & 1& &1&1& &  3& 2& 1\cr
 \phi_{80,13} &      6& & &1&1& & & & & & &  &1& & &1&  3& 2& 2\cr
 \phi_{10,14} &  6\ird& & & & &1& & & & & & 1& & & & &  1&  &  \cr
 \phi_{20,29}'&     -6& & & & & &1& & & & &  &1& & & &  1& 1&  \cr
 \phi_{20,20} & -6\ird& & & & & & &1& & & &  & &1& & &   & 1& 1\cr
 \phi_{30,20}'&  6\ird& & & & &1&1& & & & &  & & &1& &  2& 1&  \cr
 \phi_{60,20} &  6\ird& & & & & & & &1& & & 1& &1&1& &  3& 2& 1\cr
 \phi_{80,17} &      6& & & & & & &1&1& & &  &1& & &1&  3& 2& 2\cr
 \phi_{30,12}'& -6\ird& & & & & & & & &1& & 1& &1& & &  1& 1& 1\cr
 \phi_{30,24} &  6\ird& & & & & & & & & &1& 1& &1& & &  1& 1& 1\cr
\phi_{60,15}''&      6& & & & & & & & &1&1&  & & & &1&  2& 2& 2\cr
 \phi_{40,14} & -3\ird& & & & & & & & & & & 1& & &1& &  3& 1&  \cr
 \phi_{40,22} &  3\ird& & & & & & & & & & & 1& & &1& &  3& 1&  \cr
\noalign{\hrule}
 ??           &  & & & & & & & & & & & *&*&*&*& & & & \cr
\noalign{\hrule}
      }}
$$

$$\vbox{\offinterlineskip\halign{\vrule height9pt depth2.5pt width0pt
         $#$\hfil\ \vrule\  & \hfil$#$\ \vrule\  & $#$\hfil\quad &
         $#$\hfil\quad& $#$\hfil\quad&
         $#$\hfil\quad& $#$\hfil\quad \vrule\ & $#$\hfil\quad & $#$\hfil\quad&
         $#$\hfil\quad\vrule& \ $#$\hfil\quad&
         $#$\hfil\quad\vrule \cr
 \chi& \hfil f_\chi\hfil & \multispan4\hfil $p=2$\hfil& &
      \multispan2\hfil $p=3$\hfil& & \hfil c \hfil& \cr
\noalign{\hrule}
 \phi_{1,40}& -6\zd^2& 1& & & & & 1& & & 1& \cr
 \phi_{1,80}&   -6\zd&  &1& & & & 1& & & 1& \cr
 \phi_{4,41}& -3-\ird&  & &1& & &  &1& & 4& \cr
 \phi_{4,61}& -3+\ird&  & & &1& &  &1& & 4& \cr
 \phi_{5,44}&   -6\zd&  &1&1& & &  & &1& 5& \cr
 \phi_{5,52}& -6\zd^2& 1& & &1& &  & &1& 5& \cr
 \phi_{6,48}&       3&  & & & &1& 1& &1& 6& \cr
\noalign{\hrule}
      }}$$

$$ \vbox{\offinterlineskip\halign{\vrule height9pt depth3pt width0pt
         $#$\hfil\ \vrule\  & \hfil$#$\ \vrule\  & $#$\hfil\  &
         $#$\hfil\  & $#$\hfil\  & $#$\hfil\  & $#$\hfil\  &
         $#$\hfil\  & $#$\hfil\  & $#$\hfil\  &
         $#$\hfil\ \vrule\ & $#$\hfil\  & $#$\hfil\  &
         $#$\hfil\  & $#$\hfil\  & $#$\hfil &
         $#$\hfil\vrule\  &
         $#$\hfil\  & $#$\hfil\  & $#$\hfil\  &
         $#$\hfil\ \vrule\  &
         \ $#$\hfil\  &  \ $#$\hfil\ \vrule \cr
 \chi& f_\chi& \multispan{8} \hfil $p=2$\hfil& & \multispan{5}
   \hfil $p=3$\hfil& & \multispan{3}
   \hfil $p=5$\hfil& & c& \cr
\noalign{\hrule}
\phi_{4,21}  &      30&1& & & & & & & & & 1& & & & & & 1& & & &  1&  \cr
\phi_{4,51}  &      30& &1& & & & & & & & 1& & & & & &  &1& & &  1&  \cr
\phi_{24,16} &  6\zd^2& & &1& & & & & & &  &1& & & & &  & & & &  1& 1\cr
\phi_{24,26} &    6\zd& & & &1& & & & & &  &1& & & & &  & & & &  1& 1\cr
\phi_{36,15} &     -30&1& & & &1& & & & & 2& &1& & & & 1& &1& &  4& 1\cr
\phi_{36,17} &   -6\zd& & & & & &1& & & & 2& &1& & & &  & & & &  4& 1\cr
\phi_{36,25} & -6\zd^2& & & & & & &1& & & 2& &1& & & &  & & & &  4& 1\cr
\phi_{36,27} &     -30& &1& & & & & &1& & 2& &1& & & &  &1& &1&  4& 1\cr
\phi_{40,24} &       6&1&1& & & & & &1& &  & & &1& & &  & & & &  5& 1\cr
\phi_{40,18} &       6&1&1& & &1& & & & &  & & &1& & &  & & & &  5& 1\cr
\phi_{60,17} &   -6\zd& & & &1& &1& & & & 1& &2& & & &  & & & &  5& 2\cr
\phi_{60,19} & -6\zd^2& & &1& & & &1& & & 1& &2& & & &  & & & &  5& 2\cr
\phi_{64,18} &       6& & & & & & & & &1&  &1& &1& & &  & & & &  6& 2\cr
\phi_{64,21} &      30& & & & & & & & &1& 2& &2& & & &  & &1&1&  6& 2\cr
\noalign{\hrule}
 ??          & &*&*& & & & & & & & *& &*& & & &  & & &  & & \cr
\noalign{\hrule}
      }}
$$

Nous avons utilis\'e ici la matrice de d\'ecomposition (pour $p=3$)
du groupe de r\'eflexion
$G_{32}\cong 3\times2\cdot U_4(2)$ pour obtenir l'information
manquante sur les quatre derni\`eres familles, via le lemme~\ref{props}(c).

\subsection{$W=G_{33}$}

Pour $W=G_{33}$, on note $\chi=\phi_{10,8}'$ tel que
$R_\chi\equiv x^{28}+x^{26}\pmod{x^{25}}$, et $\chi=\phi_{40,5}'$ tel que
$R_\chi\equiv x^{31}+x^{29}\pmod{x^{28}}$, et on note
$\phi_{10,17}'=\phi_{10,8}\otimes\eps$, $\phi_{40,14}'=\phi_{40,5}\otimes\eps$,
avec $\eps=\phi_{1,45}$. 

Le groupe de r\'eflexions $G_{33}$ a~19 familles de caract\`eres, des types
suivants:

$\bullet$\ d\'efaut $0$: $(\phi_{1,0}),(\phi_{15,2}),(\phi_{81,6}),
 (\phi_{60,7}),(\phi_{15,9}),(\phi_{60,10}),(\phi_{81,11}),(\phi_{15,12}),
 (\phi_{15,23})$,\break $(\phi_{1,45})$

$\bullet$ type $\alpha_{1^2}(3)$~: $(\phi_{5,3},\phi_{5,1})$, $(\phi_{45,9},
  \phi_{45,7})$,
 $(\phi_{45,10},\phi_{45,12})$, $(\phi_{5,28},\phi_{5,30})$ 

$\bullet$ type $\alpha_{1^2}(2^2)$~: $(\phi_{64,8},\phi_{64,9})$

$\bullet$ type $\cM(Z_2)$~: $(\phi_{6,5},\phi_{24,4},\phi_{30,3})$,
 $(\phi_{6,20},\phi_{24,19},\phi_{30,18})$

$\bullet$ type $\alpha_{3\times5}$~: 
 $(\phi_{10,8}',\phi_{10,8}'',\phi_{20,6},\phi_{30,6},\phi_{30,4},
  \phi_{40,5}',\phi_{40,5}'')$,
 $(\phi_{10,17}'',\phi_{10,17}',\phi_{20,15},\phi_{30,13},$\\
 $\phi_{30,15},\phi_{40,14}'',\phi_{40,14}')$

\subsection{$W=G_{34}$}

Les notations utilis\'ees pour les caract\`eres irr\'eductibles de
$G_{34}$ sont celles de \cite[4C]{MaG}.

Le groupe de r\'eflexions $G_{34}$ a~40 familles de caract\`eres, des types
suivants:

$\bullet$\ d\'efaut $0$: $(\phi_{1,0}),(\phi_{70,9}'),(\phi_{210,12}),
  (\phi_{210,30}),(\phi_{70,45}'),(\phi_{1,126})$

$\bullet$ type $\alpha_{1^2}(3)$: $(\phi_{6,5},\phi_{6,1})$,
  $(\phi_{21,4},\phi_{21,2})$, $(\phi_{630,13},\phi_{630,11})$,
  $(\phi_{630,14},\phi_{630,16})$, $(\phi_{630,22},\!$ $\phi_{630,20})$,
  $(\phi_{630,23},\phi_{630,25})$, $(\phi_{21,68},\phi_{21,70})$,
  $(\phi_{6,85},\phi_{6,89})$

$\bullet$ type $\cM(Z_2)$~: $(\phi_{21,6},\phi_{35,6},\phi_{56,3})$,
  $(\phi_{140,12},\phi_{420,12},\phi_{560,9})$,
  $(\phi_{140,21},\phi_{420,21},\phi_{560,18}'),$\\
  $(\phi_{140,30},\phi_{420,30},\phi_{560,27})$,
  $(\phi_{21,60},\phi_{35,60},\phi_{56,57})$

$\bullet$ type $\alpha_{1^2}(2^2)$~: $(\phi_{896,12},\phi_{896,15})$,
  $(\phi_{896,21},\phi_{896,24})$

$\bullet$ type $\alpha_{3\times5}$~: $(\phi_{15,16},\phi_{15,14},\phi_{90,6},
  \phi_{105,8}',\phi_{105,4},\phi_{120,7},\phi_{120,5})$,
 $(\phi_{15,56},\phi_{15,58},\phi_{90,48}$,\\
 $\phi_{105,46},\phi_{105,50},\phi_{120,47},\phi_{120,49})$

$\bullet$ type $\cM(Z_3)$~: $(\phi_{70,9}''',\phi_{70,9}'',\phi_{56,9},
  \phi_{126,5},\phi_{126,7})$,
 $(\phi_{280,12}'',\phi_{280,12}',\phi_{35,18},\phi_{315,10},\phi_{315,14}),\!\!$\\
 $(\phi_{210,13},\phi_{210,17},\phi_{630,15},\phi_{840,11},\phi_{840,13}')$,
 $(\phi_{210,29},\phi_{210,25},\phi_{630,27},\phi_{840,23}'',\phi_{840,25})$,
 \\
 $(\phi_{280,30}'$, $\phi_{280,30}'',\phi_{35,36},\phi_{315,28},\phi_{315,32})$,
 $(\phi_{70,45}'',\phi_{70,45}''',\phi_{56,45},\phi_{126,41},\phi_{126,43})$

$\bullet$ type $\cM(Z_3)''$~: $(\phi_{105,8}',\phi_{105,10},\phi_{210,8},
  \phi_{210,10},\phi_{315,6})$,
 $(\phi_{105,38},\phi_{105,40},\phi_{210,38},\phi_{210,40},\ \ $
 $\phi_{315,36})$

$\bullet$ type $\alpha_{3\times5}''$~: $(\phi_{84,13},\phi_{84,17},\phi_{336,8},
  \phi_{336,10},\phi_{420,11},\phi_{420,7},\phi_{504,9})$,
 $(\phi_{84,41},\phi_{84,37},\phi_{336,34},\ $ 
 $\phi_{336,32},\phi_{420,35},\phi_{420,33},\phi_{504,31})$,

$\bullet$ type $\alpha_{3\times4}(3)$~: $(\phi_{384,8},\phi_{384,11},
 \phi_{384,13}, \phi_{384,10})$, $(\phi_{384,34},\phi_{384,31},\phi_{384,29},
 \phi_{384,32})$

$\bullet$ type $\alpha_{1^3}(3)$~: $(\phi_{729,10},\phi_{729,12},\phi_{729,14})$,
 $(\phi_{729,24},\phi_{729,26},\phi_{729,28})$,

$\bullet$ Autres familles~:

$$ \vbox{\offinterlineskip\halign{\vrule height9pt depth3pt width0pt
         $#$\hfil\ \vrule\ & $#$\hfil\ \vrule\ & \hfil$#$\ \vrule\ &
         $#$\hfil\  &
         $#$\hfil\  & $#$\hfil\  & $#$\hfil\  & $#$\hfil\  &
         $#$\hfil\  & $#$\hfil\  & $#$\hfil\  & $#$\hfil\  &
         $#$\hfil\ \vrule\  & $#$\hfil\  & $#$\hfil\  &
         $#$\hfil\  & $#$\hfil\  & $#$\hfil\  &
         $#$\hfil\ \vrule\ &
         \ $#$\hfil\  & \ $#$\hfil\  & \ $#$\hfil\  &
         \ $#$\hfil\ \vrule\cr
 \chi& \bar\chi\otimes\eps& f_\chi& \multispan{9} \hfil $p=2$\hfil& & \multispan{5}
   \hfil $p=3$\hfil& & & c& & \cr
\noalign{\hrule}
 \phi_{105,20}& \phi_{105,26} &      -6& 1& & & & & & & & & & 1& & & & & &  1&  &  &  \cr
 \phi_{105,22}& \phi_{105,28} &      -6&  &1& & & & & & & & & 1& & & & & &  1&  &  &  \cr
 \phi_{189,18}& \phi_{189,24} &       6&  & &1& & & & & & & &  &1& & & & &   & 1&  &  \cr
 \phi_{315,18}& \phi_{315,24} &       6&  & & &1& & & & & & &  & &1& & & &   &  & 1&  \cr
 \phi_{336,17}& \phi_{336,25} & -6\zd^2&  & & & &1& & & & & &  & & &1& & &  1&  &  & 1\cr
 \phi_{336,19}& \phi_{336,23} &   -6\zd&  & & & & &1& & & & &  & & &1& & &  1&  &  & 1\cr
 \phi_{420,14}& \phi_{420,22} &  6\zd^2&  & & & & & &1& & & & 1& &1& & & &  1&  & 1&  \cr
 \phi_{420,16}& \phi_{420,20} &    6\zd&  & & & & & & &1& & & 1& &1& & & &  1&  & 1&  \cr
 \phi_{504,15}& \phi_{504,21} &       6&  & &1&1& & & & & & &  & & & &1& &   & 1& 1&  \cr
 \phi_{756,14}& \phi_{756,22} & -6\zd^2&  & & & &1& &1& & & &  & & & & &1&  2&  & 1& 1\cr
 \phi_{756,16}& \phi_{756,20} &   -6\zd&  & & & & &1& &1& & &  & & & & &1&  2&  & 1& 1\cr
 \phi_{840,13}''& \phi_{840,19}&      6&  & & & & & & & &1& &  & & &1&1& &  1& 1& 1& 1\cr
 \phi_{840,17}& \phi_{840,23}'&       6&  & & & & & & & & &1&  & & &1&1& &  1& 1& 1& 1\cr
 \phi_{945,14}& \phi_{945,20} &       6& 1& & & & & & & & &1&  &1& & & &1&  2& 1& 1& 1\cr
 \phi_{945,16}& \phi_{945,22} &       6&  &1& & & & & & &1& &  &1& & & &1&  2& 1& 1& 1\cr
\noalign{\hrule}
      }}
$$

$$ \vbox{\offinterlineskip\halign{\vrule height9pt depth3pt width0pt
         $#$\hfil\ \vrule\ & \hfil$#$\ \vrule\ &  $#$\hfil\  &
         $#$\hfil\  & $#$\hfil\  & $#$\hfil\  &
         $#$\hfil\  & $#$\hfil\  & $#$\hfil\  &
         $#$\hfil\  & $#$\hfil\ \vrule\ & 
         $#$\hfil\  & $#$\hfil\  & $#$\hfil\  &
         $#$\hfil\  & $#$\hfil\ \vrule\ & 
         $#$\hfil\  & $#$\hfil\  & $#$\hfil\  &
         $#$\hfil\  & $#$\hfil\ \vrule\ &
         \ $#$\hfil\  & \ $#$\hfil\  & \ $#$\hfil\  &
         $#$\hfil\ \vrule\cr
 \chi& f_\chi& \multispan{8} \hfil $p=2$\hfil& & \multispan{4}
   \hfil $p=3$\hfil& & \multispan{4} \hfil $p=7$\hfil& & & c& & \cr
\noalign{\hrule}
\phi_{20,33}'   &  42&1& & & & & & & & & 1& & & & & 1& & & & &  1&  &  & \cr
\phi_{20,33}''  &  42& &1& & & & & & & & 1& & & & &  &1& & & &  1&  &  & \cr
\phi_{120,21}'  &   7& & & & & & & & & &  & & & & &  & &1& & &   & 1&  & \cr
\phi_{120,21}'' &   7& & & & & & & & & &  & & & & &  & &1& & &   & 1&  & \cr
\phi_{540,17}   &   6& & &1& & & & & & & 2&1& & & &  & & & & &  4&  & 1& \cr
\phi_{540,19}   &   6& & & &1& & & & & & 2&1& & & &  & & & & &  4&  & 1& \cr
\phi_{540,21}'  & -42& &2& & &1& & & & & 2&1& & & &  &1& &1& &  4&  & 1& \cr
\phi_{540,21}'' & -42&2& & & & &1& & & & 2&1& & & & 1& & & &1&  4&  & 1& \cr
\phi_{560,18}'' &   6&1&2& & &1& & & & &  & &1& & &  & & & & &  5&  & 1& \cr
\phi_{560,18}'''&   6&2&1& & & &1& & & &  & &1& & &  & & & & &  5&  & 1& \cr
\phi_{720,16}   &   6& & & & & & &1& & &  & & &1& &  & & & & &  1& 2& 1& \cr
\phi_{720,20}   &   6& & & & & & & &1& &  & & &1& &  & & & & &  1& 2& 1& \cr
\phi_{1260,17}  &   6& & &1& & & & &1& & 1&2& & &1&  & & & & &  5& 2& 2& \cr
\phi_{1260,19}  &   6& & & &1& & &1& & & 1&2& & &1&  & & & & &  5& 2& 2& \cr
\phi_{1280,18}  &   6& & & & & & & & &1&  & &1&1& &  & & & & &  6& 2& 2& \cr
\phi_{1280,15}  &  42& & & & & & & & &1& 2&2& & &1&  & &2&1&1&  6& 2& 2& \cr
\noalign{\hrule}
??              &    & & & & & & & & & & *&*& & & &  & & & & &  & & & \cr
\noalign{\hrule}
      }}
$$

Pour justifier ces r\'esultats, il nous faut d\'emontrer
que la r\'eduction modulo $3$ de $\phi_{420,11}$ reste simple.
Soit $C$ le bloc de $\cO_\fp\cH$ contenant ce caract\`ere.

Soit $W'$ un sous-groupe parabolique de type $G_{33}$.
Soit $L$ un $k_\fp\cO(W')$ module simple, de caract\`ere 
$d(\phi_{10,8}')$. Soit $C'$ le bloc de
$\cO_\fp\cH(W')$ contenant $\phi_{10,8}'$. La projection sur $C$ de 
$\Ind_{W'}^W \phi_{10,8}'$ est $\phi_{84,17}$, dont la r\'eduction
modulo $3$ est un module simple $S$. Par cons\'equent,
$C\otimes_{C'}L\simeq S$. Soit $T$ l'autre $C$-module simple.
On a $0=\Hom_C(S,T)\simeq \Hom_{C'}(L,C'\otimes_C T)$.
Par cons\'equent, $C'\otimes_C T$ n'a pas de sous-module isomorphe \`a
$S$. On montre de m\^eme qu'il n'a pas de quotient isomorphe \`a 
$S$.
Soit $X$ un $C$-module, libre
sur $\cO_\fp$, de caract\`ere $\phi_{420,11}$. Supposons sa r\'eduction
modulo $\fp$ non simple~: elle a alors deux facteurs de composition $S$
et $T$. Alors, $C'\otimes_C T$ a deux facteurs de composition distincts, $L$
et $L'$. Mais un tel module a n\'ecessairement un sous-module ou
un quotient isomorphe \`a $L$, ce qui est impossible.

On montre de m\^eme que la reduction modulo $3$ de
$\phi_{420,35}$ est simple.

\subsection{$W=G_{35}=W(E_6)$}

Le groupe de r\'eflexions $G_{35}$ a~17 familles de caract\`eres, des types
suivants (voir \cite{Gy}):

$\bullet$ type $\cM(Z_2)$~: $(\phi_{30,3},\phi_{15,5},\phi_{15,4})$,
  $(\phi_{30,15},\phi_{15,17},\phi_{15,16})$.

$\bullet$ type $\cM(S_3)$~: $(\phi_{10,9},\phi_{20,10},\phi_{60,8},\phi_{80,7},
  \phi_{90,8})$

Tous les autres familles sont de d\'efaut $0$.

\subsection{$W=G_{36}=W(E_7)$}

Le groupe de r\'eflexions $G_{36}$ a~35 familles de caract\`eres, des types
suivants (voir \cite{Gy}):

$\bullet$ type $\alpha_{1^2}(2^2)$~: $(\phi_{512,11},\phi_{512,12})$

$\bullet$ type $\cM(Z_2)$~: $(\phi_{56,3},\phi_{35,4},\phi_{21,6})$,
$(\phi_{405,8},\phi_{216,9},\phi_{189,10})$,
$(\phi_{420,10},\phi_{84,12},\phi_{336,11})$,\\
$(\phi_{420,13},\phi_{84,15},\phi_{336,14})$,
$(\phi_{405,15},\phi_{216,16},\phi_{189,17})$,
$(\phi_{56,30},\phi_{35,31},\phi_{21,33})$

$\bullet$ type $\cM(S_3)$~: $(\phi_{70,9},\phi_{35,13},\phi_{280,8},
\phi_{315,7},\phi_{280,9})$,
$(\phi_{70,18},\phi_{35,22},\phi_{280,17},\phi_{315,16},\phi_{280,18})$

Toutes les autres familles sont de d\'efaut $0$.

\subsection{$W=G_{37}=W(E_8)$}

Le groupe de r\'eflexions $G_{37}$ a~46 familles de caract\`eres.
Gyoja a montr\'e que 23 d'entre elles correspondent \`a des blocs de
d\'efaut~0, seize sont de type $\cM(Z_2)$.

$\bullet$ type $\alpha_{1^2}(2^2)$~: $(\phi_{4096,11},\phi_{4096,12})$,
$(\phi_{4096,26},\phi_{4096,27})$

$\bullet$ type $\cM(S_3)$~:
$(\phi_{448,9},\phi_{56,19},\phi_{1344,8},\phi_{1400,7},\phi_{1008,9})$,
$(\phi_{175,12},\phi_{350,14},\phi_{1050,10},\phi_{1400,8},$\\
$\phi_{1575,10})\!$,
$\!(\phi_{175,36},\phi_{350,38},\phi_{1050,34},\phi_{1400,32},\phi_{1575,34})\!$,
$\!(\phi_{448,39},\phi_{56,49},\phi_{1344,38},\phi_{1400,37},\phi_{1008,39})\!\!$

La derni\`ere famille \`a 17 \'el\'ements a la matrice de
d\'ecomposition suivante:

$$ \vbox{\offinterlineskip\halign{\vrule height9pt depth3pt width0pt
         $#$\hfil\quad\vrule\ & \hfil$#$\ \vrule\ & $#$\hfil\  & $#$\hfil\  &
         $#$\hfil\  & $#$\hfil\  & $#$\hfil\  & $#$\hfil\  &
         $#$\hfil\  & $#$\hfil\  & $#$\hfil\  & $#$\hfil\ \vrule\ & 
         $#$\hfil\  & $#$\hfil\  & $#$\hfil\  & $#$\hfil\  &
         $#$\hfil\  & $#$\hfil\   & $#$\hfil\  &
         $#$\hfil\ \vrule\ & $#$\hfil\  & $#$\hfil\  & $#$\hfil\  &
         $#$\hfil\ \vrule\cr
\phi& f_\chi&  \multispan{9} \hfil $p=2$\hfil& & \multispan{7}
   \hfil $p=3$\hfil& & \multispan{3} \hfil $p=5$\hfil& \cr
\noalign{\hrule}
 \phi_{70,32}&   30&1& & & & & & & & & & 1& & & & & & & & 1& & & \cr
 \phi_{168,24}&   8& &1& & & & & & & & &  & & & & & & & &  & & & \cr
 \phi_{420,20}&   5& & &1& & & & & & & &  & & & & & & & &  &1& & \cr
 \phi_{448,25}&  12& & & &1& & & & & & &  &1& & & & & & &  & & & \cr
 \phi_{1134,20}&  6& & & & &1& & & & & &  & &1& & & & & &  & & & \cr
 \phi_{1344,19}&  4& & & & & &1& & & & &  & & & & & & & &  & & & \cr
 \phi_{1400,20}& 24& &1& & & & &1& & & & 1& & &1& & & & &  & & & \cr
 \phi_{1680,22}& 20& & & &1& & &1& & & &  & & & & & & & & 1& &1& \cr
 \phi_{2016,19}&  6& & & & & & & &1& & &  & & & &1& & & &  & & & \cr
 \phi_{2688,20}&  8& & & & & & & & &1& &  & & & & & & & &  & & & \cr
 \phi_{3150,18}&  6& & & & &1& & &1& & &  & & & & &1& & &  & & & \cr
 \phi_{4200,18}&  8& &1& & & &1& & &1& &  & & & & & & & &  & & & \cr
 \phi_{4480,16}&120& & & &1& &1& & &1& &  & & &1& &1& & &  &1& &1\cr
 \phi_{4536,18}& 24& &1& &1& & &1& &1& &  & & & & & &1& &  & & & \cr
 \phi_{5600,19}&  6& & & & & & & & & &1&  &1& & & & & &1&  & & & \cr
 \phi_{5670,18}& 30&1& & & & & & & & &1&  & &1& & & &1& &  & &1&1\cr
 \phi_{7168,17}& 12& & & &1& &1& & &2& &  & & & &1& & &1&  & & & \cr
\noalign{\hrule}
      }}
$$

Pour les sept caract\`eres constructibles, voir Lusztig \cite{Lu1}.

\subsection{Groupes di\'edraux}

Dans cette partie, nous d\'eterminons les blocs pour les
groupes di\'edraux $I_2(n)$,
$n\ge3$. Les caract\`eres des groupes di\'edraux se r\'epartissent en trois
familles: deux familles \`a un \'el\'ement, consistant en le
caract\`ere trivial et le caract\`ere signe, et une famille contenant
les autres caract\`eres. Les \'el\'ements de Schur pour les groupes
de r\'eflexions di\'edraux ont \'et\'e d\'etermin\'es par
Kilmoyer et Solomon (cf. par exemple \cite[Thm.~8.3.4]{GP}).
Ce r\'esultat a aussi \'et\'e obtenu par Brou\'e et Kim \cite{BrKi}
et il r\'esulte aussi de \cite{Ro} dont les m\'ethodes restent
valables pour les groupes di\'edraux (car on dispose encore de
l'alg\`ebre $J$ de Lusztig). La preuve que nous donnons ici
est beaucoup plus directe.

\begin{thm} \label{diedral}
 Soit $W=I_2(n)$, $n\ge3$. Alors,
 les blocs de  $\cO\cH(W)$ co\"{\i}ncident avec les familles de Lusztig.
\end{thm}

\begin{proof}[Preuve]
Pour $\chi$ le caract\`ere trivial ou le caract\`ere signe, on a
$c_\chi\in\cO^\times$. D'apr\`es le lemme \ref{props}(b), 
il suffit donc de prouver que les autres caract\`eres sont dans un
m\^eme bloc. \par

Supposons tout d'abord que $n=2m+1$ avec $m$ entier.
Soit $\zeta$ une racine primitive $n$-\`eme de l'unit\'e.
Hormis les caract\`eres triviaux et signes, il y a
$m$ caract\`eres $\rho_i$, $1\le i\le m$, tous de degr\'e $2$.
On a
$f_{\rho_i}^{-1}=a_i$, o\`u
$$a_i=\frac{(1-\zeta^i)(1-\zeta^{-i})}{n}.$$
Par cons\'equent, $|a_i|<1$ et $\sum_{i=1}^m a_i =1$.
L'induite de la repr\'esentation triviale \`a partir
d'un sous-groupe parabolique d'ordre $2$ fournit un caract\`ere constructible
dans lequel tous les $\rho_i$ apparaissent avec multiplicit\'e $1$.
Soit $I$ un sous-ensemble propre non vide de $\{1,\ldots,m\}$ et
$a_I=\sum_{i\in I} a_i$. Alors, $|a_I|<1$ et le module d'un conjugu\'e
galoisien de $a_I$ est aussi strictement inf\'erieur \`a $1$ (c'est encore
une somme partielle de $a_i$). En particulier, la $\QQ(\zeta)/\QQ$-norme
de $a_I$ n'est pas enti\`ere, donc $a_I$ n'est pas un entier.
Par cons\'equent, $\sum_i \rho_i$ est
un caract\`ere projectif ind\'ecomposable, d'apr\`es le lemme~\ref{props}(a).

Supposons maintenant $n=2m$ avec $m$ entier et soit $\zeta$ une
racine primitive $n$-\`eme de l'unit\'e.
Les caract\`eres diff\'erents des caract\`eres triviaux et signes consistent
en deux caract\`eres lin\'eaires (avec $f_{\chi}^{-1}=2/n$) et
$m-1$ caract\`eres
$\rho_i$, $1\le i\le m-1$, de degr\'e $2$ (avec $f_{\rho_i}^{-1}=a_i$).

En induisant les deux caract\`eres d'un sous-groupe parabolique d'ordre
$2$, on obtient deux caract\`eres projectifs contenant chacun
l'un des deux caract\`eres de degr\'e $1$ et contenant les $\rho_i$ avec
multiplicit\'e $1$. Comme pr\'ec\'edemment, on voit que ce sont
des caract\`eres projectifs ind\'ecomposables.
On en d\'eduit la co\"{\i}ncidence entre blocs de $\cO\cH$ et familles.
\end{proof}

Dans \cite{Ro}, ce r\'esultat est \'etabli pour les groupes de Weyl 
ou plus g\'en\'eralement lorsqu'on dispose de l'alg\`ebre $J$ de Lusztig,
ce qui est connu sauf pour le type $H_4$.

\begin{cor} \label{family=block}
 Soit $W$ un groupe de Coxeter fini. Alors, les blocs de l'alg\`ebre de
 Hecke de $W$ sur $\cO$ co\"{\i}ncident avec les familles de Lusztig.
\end{cor}

\section{Quelques groupes non spetsiaux}

Dans cette partie, nous d\'eterminons les matrices de d\'ecomposition
et la division en familles pour quelques groupes de r\'eflexions complexes
exceptionnels non spetsiaux. Nous utiliserons plus tard ces r\'esultats
pour montrer que certaines propri\'et\'es des familles des groupes spetsiaux
ne s'\'etendent pas au cas non spetsial. Les techniques d\'evelopp\'ees
au \S \ref{sec:explicit} s'appliquent sans changement au cas non spetsial. Il
serait int\'eressant de d\'eterminer les matrices de d\'ecomposition
pour tous les groupes non spetsiaux et tous les choix de param\`etres.

\subsection{$W=G_5$} \label{subsec:G5}

Le groupe de r\'eflexions complexe 
$G_5$ est d'ordre 72 et poss\`ede 21 caract\`eres irr\'educibles.
Les matrices de d\'ecomposition de $\cO_\fp\cH(W)$ pour les mauvais nombres
premiers 2 et 3 et la division en familles sont les suivantes:

$\bullet$\ d\'efaut $0$~: $(\phi_{1,0})$

$\bullet$ type $\alpha_{1^2}(3)$~:
$(\phi_{2,5}',\phi_{2,7}')$, $(\phi_{2,5}''',\phi_{2,7}'')$

$\bullet$ Autres familles~:
$$ \vbox{\offinterlineskip\halign{\vrule height9pt depth2.5pt width0pt
         $#$\hfil\quad\vrule\  & $#$\hfil\quad\vrule\ & $#$\hfil\quad\vrule\
         & $#$\hfil\quad\vrule \cr
  &\hfil f_\chi\hfil & \hfil p=3 \hfil & \hfil c \hfil \cr
\noalign{\hrule}
\phi_{3,6} &      3& 1& 1\cr
\phi_{3,2} & 3\zd  & 1& 1\cr
\phi_{3,4} & 3\zd^2& 1& 1\cr
\noalign{\hrule}
      }}
$$

$$ \vbox{\offinterlineskip\halign{\vrule height9pt depth2.5pt width0pt
         $#$\hfil\ \vrule\  & \hfil$#$\ \vrule\ & $#$\hfil\quad &
         $#$\hfil\quad & $#$\hfil\quad & $#$\hfil\quad\vrule\  &
         $#$\hfil\quad & $#$\hfil\quad\vrule& \ $#$\hfil\quad\vrule \cr
 &\hfil f_\chi\hfil & \multispan3 \hfil $p=2$\hfil& &
    \multispan{2}\hfil $p=3$ \hfil\vrule& c \cr
\noalign{\hrule}
\phi_{1,8}''  &3\zd^2 & 1&  &  &  &  1&  & 1\cr
\phi_{1,12}'  &6      &  & 1&  &  &  1&  & 1 \cr
\phi_{1,12}'' &6      &  &  & 1&  &  1&  & 1\cr
\phi_{1,16}   &3\zd   &  &  &  & 1&  1&  & 1\cr
\phi_{2,9}    &2      &  & 1& 1&  &   & 1& 2\cr
\noalign{\hrule}
      }}
$$

$$ \vbox{\offinterlineskip\halign{\vrule height9pt depth2.5pt width0pt
         $#$\hfil\quad\vrule\  & \hfil$#$\quad\vrule\  &
         $#$\hfil\quad & $#$\hfil\quad &
         $#$\hfil\quad & $#$\hfil\quad &
         $#$\hfil\quad & $#$\hfil\quad\vrule\ &
         $#$\hfil\quad & $#$\hfil\quad\vrule\ & $#$\hfil\quad &
         $#$\hfil\quad\vrule \cr
 & f & \multispan5 \hfil $p=2$\hfil& & \multispan1\hfil $p=3$\hfil& &
 c&  \cr
\noalign{\hrule}
\phi_{1,4}'   &2(\zd^2-\zd) & 1 &   &   &   &   &  &  1&   & 1&  \cr
\phi_{1,4}''  &2(\zd^2-\zd) &   & 1 &   &   &   &  &   & 1 &  & 1\cr
\phi_{1,8}'   &2(\zd-\zd^2) &   &   & 1 &   &   &  &  1&   & 1&  \cr
\phi_{1,8}''' &2(\zd-\zd^2) &   &   &   & 1 &   &  &   & 1 &  & 1\cr
\phi_{2,5}''  &6            &   &   & 1 & 1 &   &  &  1& 1 & 1& 1\cr
\phi_{2,3}'   &3            &   &   &   &   & 1 &  &  1& 1 & 1& 1\cr
\phi_{2,3}''  &3            &   &   &   &   &   & 1&  1& 1 & 1& 1\cr
\phi_{2,1}    &6            & 1 & 1 &   &   &   &  &  1& 1 & 1& 1\cr
\noalign{\hrule}
      }}
$$

\subsection{$W=G_{12}$}

Le groupe de r\'eflexions complexe non spetsial
$G_{12}$ a les cinq familles suivantes:

$\bullet$\ d\'efaut $0$~: $(\phi_{1,0}),(\phi_{3,2}),(\phi_{3,6}),
(\phi_{1,12})$

$\bullet$ Autre famille~:
$$ \vbox{\offinterlineskip\halign{\vrule height9pt depth2.5pt width0pt
         $#$\hfil\quad\vrule\ & \hfil$#$\quad\vrule\ &  $#$\hfil\ &
         $#$\hfil\quad\vrule\ & $#$\hfil\quad & $#$\hfil\quad\vrule &
         \ $#$\hfil\quad\vrule \cr
 &  f_\chi& \multispan1\hfil $p=2$\hfil& & \multispan2\hfil $p=3$\hfil\vrule&
   \hfil c\cr
\noalign{\hrule}
\phi_{2,1} & -12& 1& &  1&  & 1\cr
\phi_{2,4} &   2& 1& &   & 1& 1\cr
\phi_{2,5} & -12& 1& &  1&  & 1\cr
\phi_{4,3} &   3&  &1&  2&  & 2\cr
\noalign{\hrule}
 ??        &    &  & &  *&  & \cr
\noalign{\hrule}
      }}
$$

\subsection{$W=G_{13}$}

Le groupe de r\'eflexions complexe non spetsial
$G_{13}$ a les sept familles suivantes:

$\bullet$\ d\'efaut $0$~: $(\phi_{1,0}),(\phi_{3,4}),(\phi_{3,6}),
 (\phi_{1,18})$

$\bullet$ type $\cM(Z_2)$~:
 $(\phi_{1,6},\phi_{2,4},\phi_{3,2})$, $(\phi_{1,12},\phi_{2,10},\phi_{3,8})$

$\bullet$ Autre famille~:
$$ \vbox{\offinterlineskip\halign{\vrule height9pt depth2.5pt width0pt
         $#$\hfil\quad\vrule\ & \hfil$#$\quad\vrule\ & $#$\hfil\quad&
         $#$\hfil\quad\vrule\ & $#$\hfil\quad & $#$\hfil\quad\vrule &
         \ $#$\hfil\quad\vrule \cr
 &  f_\chi& \hfil p=2\hfil& & \multispan2\hfil $p=3$\hfil\vrule& c\cr
\noalign{\hrule}
\phi_{2,7}' & -12\zeta_4& 1& &  1&  & 1\cr
\phi_{2,1}  & 12\zeta_4&  1& &   & 1& 1\cr
\phi_{2,5}  & 12\zeta_4&  1& &   & 1& 1\cr
\phi_{2,7}''& -12\zeta_4& 1& &  1&  & 1\cr
\phi_{4,3}  & -6\zeta_4&  2& &  2&  & 2\cr
\phi_{4,5}  &  6\zeta_4&  2& &   & 2& 2\cr
\noalign{\hrule}
 ??         &           & *& &  *& *& \cr
\noalign{\hrule}
      }}
$$

\section{Quelques cons\'equences}

Les r\'esultats de la partie~\ref{sec:explicit}, joints au r\'esultat du
second auteur \cite{Ro} pour les groupes de Weyl et aux travaux de Brou\'e et
Kim \cite{BrKi} sur les groupes imprimitifs
(s\'eries infinies $G(d,1,n)$ et $G(d,d,n)$) compl\`etent la
d\'etermination des familles de tous les groupes spetsiaux.

Nous d\'eduisons de cette description des propri\'et\'es g\'en\'erales,
satisfaites par toutes les familles de groupes de r\'eflexions spetsiaux
--- nous ne poss\`edons pour le moment aucune preuve directe de ces
propri\'et\'es.

Soit $W$ un groupe de r\'eflexions complexe, de corps de
d\'efinition un corps de nombre $k$.
Soit $\mu(k)$ l'ordre du groupe des racines de l'unit\'e de
$k$ et $y$ tel que $y^{\mu(k)}=x$. Alors, sous l'hypoth\`ese que la
conjecture~\ref{conj1} est vraie, le r\'esultat principal de
\cite{MaR} montre que l'alg\`ebre $K\cH(W)$ est d\'eploy\'ee pour $K=k(y)$.
Par cons\'equent,
les \'el\'ements de Schur $c_\chi$ ($\chi\in\Irr(W)$) sont contenus dans $K$.
\par
Pour $\chi\in\Irr(W)$, soit $\delta_\chi=P(W)/c_\chi$ le {\em degr\'e
g\'en\'erique} de $\chi$, une fonction rationelle de $y$. Soient
$a(\chi)\in \frac{1}{\mu(k)}\ZZ$ le plus grand nombre tel
que $y^{-a(\chi)\mu(k)}\delta_\chi$ est un polyn\^ome en $y$ et
$A(\chi)\in \frac{1}{\mu(k)}\ZZ$ le plus petit nombre tel que
$y^{A(\chi)\mu(k)}\delta_\chi(y^{-1})$ est un polyn\^ome en $y$. \par
Par d\'efinition, si de plus $W$ est spetsial, les $c_\chi$ sont m\^eme
rationnels, c'est-\`a-dire, contenus dans $k(x)$, et alors
$a(\chi),A(\chi)\in\ZZ$. \par

D'apr\`es Brou\'e--Malle--Michel \cite[Cor.~6.9]{spets}, on a
$$a(\chi)+A(\chi)= \frac{N(\chi)+N(\chi^*)}{\chi(1)},$$
donc, $a(\chi)+A(\chi)$ est constant sur les familles par le
lemme~\ref{central}. Les r\'esultats explicites pr\'ec\'edents d\'emontrent le
suivant, bien connu pour les groupes de Coxeter (voir Lusztig \cite{Lu1}), et
pour les groupes de r\'eflexions imprimitifs (voir
Brou\'e--Kim \cite[Prop.~4.5]{BrKi})~:

\begin{thm} \label{aconst}
 Soit $W$ un groupe de r\'eflexions complexe spetsial. Supposons les
 conjectures~\ref{conj1}, \ref{conj2} et \ref{conj3} vraies si
 $W$ est de type exceptionnel. Alors $a(\chi)$ (et $A(\chi)$) est constant
 sur toute famille de $W$.
\end{thm}

Notons que nous ne connaissons pas de contre-exemple au
th\'eor\`eme~\ref{aconst} pour des groupes non spetsiaux.

\smallskip
Le th\'eor\`eme \ref{aconst} ne se d\'eduit pas du
lemme~\ref{central} qui affirme la constance de $a+A$, comme le montre
l'exemple suivant:

\begin{exmp} \label{Ex1}
 {\rm L'ensembles des caract\`eres du groupe
 $G_{26}$ avec $a+A=30$ est la r\'eunion de deux familles,
 l'une avec $a(\cF)=4$ et l'autre avec $a(\cF)=5$.}
\end{exmp}

Le th\'eor\`eme~\ref{aconst} permet de d\'efinir $a(\cF)$ pour une famille
$\cF$ comme valeur commune des $a(\chi)$ pour $\chi\in\cF$. De m\^eme, nous
posons $A(\cF)=A(\chi)$ pour $\chi\in\cF$. \par
Un caract\`ere $\chi\in\Irr(W)$ est dit {\em sp\'ecial} si
$$a(\chi)=b(\chi)\,,$$
o\`u, comme not\'e dans 3.3, $b(\chi)$ est le plus grand entier tel que
$x^{-b(\chi)}R_\chi$ est un polyn\^ome. Les caract\`eres sp\'eciaux ont
\'et\'e determin\'es dans \cite[8B]{MaG}. 

\begin{thm} \label{spec}
 Soit $W$ un groupe de r\'eflexions
 complexe spetsial. Supposons les conjectures~\ref{conj1}, \ref{conj2}
 et \ref{conj3} vraies si
 $W$ est de type exceptionnel. Alors chaque famille contient un unique
 caract\`ere sp\'ecial.
\end{thm}

Pour les s\'eries infinies, cela d\'ecoule de
\cite[Bem.~2.24 et Lemma~5.16]{MaU}
et du fait que les familles combinatoires co\"\i ncident avec les blocs
\cite[Thm.~3.16]{BrKi}. Ce r\'esultat a \'et\'e v\'erif\'e (cas par cas)
pour les groupes de Coxeter par Lusztig \cite{Lu1}. \par

La proposition suivante montre que pour $W$ spetsial les
caract\`eres sp\'eciaux peuvent aussi \^etre d\'efinis comme les caract\`eres
$\chi$ tels que $b(\chi)=a(\cF(\chi))$.

D'apr\`es \cite[Prop.~8.1]{MaG}, on a:

\begin{prop} \label{aA}
 Soient $W$ un groupe de r\'eflexions complexe
 spetsial et $\cF$ une famille de $W$. Alors
 $$a(\cF)\le b(\chi)\qquad\text{et}\qquad
   B(\chi^*)\le A(\cF)$$
 pour tous $\chi\in\cF$. De plus, les in\'egalit\'es sont strictes sauf si
 $\chi$ est le caract\`ere sp\'ecial de $\cF$.
\end{prop} 

Pour les groupes non spetsiaux, les assertions~\ref{spec} et~\ref{aA}
deviennent fausses~:

\begin{exmp} \label{Ex2}
 {\rm
Les familles de $W=G_5$ sont d\'ecrites dans~\ref{subsec:G5}.
La fonction $a$ est constante sur les familles.
Pour les deux familles de type
$\alpha_{1^2}(3)$, on a $a(\cF)=4$, mais ces familles ne contiennent pas de
caract\`ere
sp\'ecial; de m\^eme pour les deux familles \`a 3 et 5 caract\`eres.
De plus, pour la famille \`a~3 caract\`eres, on a
$a(\cF)=8/3>2=b(\phi_{3,2})$. \par
De m\^eme, pour la famille $(\phi_{3,2})$ de $G_{12}$, on a $a(\cF)=1$.
}
\end{exmp}

Notons aussi que:

\begin{prop} \label{galois}
 Soit $W$ un groupe de r\'eflexions complexe spetsial. Supposons les
 conjectures~\ref{conj1}, \ref{conj2} et \ref{conj3} vraies si
 $W$ est de type exceptionnel. Alors toutes les familles de $\Irr(W)$ sont
 invariantes par les automorphismes galoisiens de $\Irr(W)$.
\end{prop} 

\begin{rem} \label{construc}
 {\rm Soit $W$ un groupe
spetsial, $\phi$ un caract\`ere
 constructible de $W$ dans la famille $\cF$ et $\chi_s$ le caract\`ere
 sp\'ecial de $\cF$. Alors, 
 $$\langle \phi,\chi_s-\sum_{\chi\in\cF} \frac{1}{f_\chi}\chi\rangle=0$$
 lorsque $W$ est exceptionnel. Cette propri\'et\'e est aussi vraie
 lorsque $W$ est un groupe de Weyl~: elle r\'esulte de l'invariance des
 caract\`eres constructibles par transformation de Fourier, elle-m\^eme
 cons\'equence de la commutation entre la transformation de Fourier et la
 $J$-induction. Il est probable que cette propri\'et\'e est vraie pour
 les groupes spetsiaux de la s\'erie infinie. Il est aussi envisageable
 que cette propri\'et\'e soit vraie plus g\'en\'eralement lorsque $\phi$
 est le caract\`ere d'un $\cO\cH$-module projectif.

 Notons aussi que lorsque
 $\sum_{\chi\in\cF} \frac{\langle \phi_1,\chi\rangle}{f_\chi}\not\in\cO$
 pour tout sous-caract\`ere propre $\phi_1$ de $\phi$, alors
 $\phi$ est le caract\`ere d'un $\cO\cH$-module projectif ind\'ecomposable.

 Dans le cas des groupes de Coxeter finis, on a $\langle \phi,\chi_s\rangle=1$. 
 Dans les cas \'etudi\'es dans~\ref{sec:explicit}, il arrive que cette
 multiplicit\'e
 soit $>1$. Ceci sugg\`ere que notre d\'efinition de caract\`eres
 constructibles n'est peut-\^etre pas la bonne dans ces cas.}
\end{rem}

\section{Compatibilit\'e entre familles, $\cO$-blocs et $d$-s\'eries de
Harish-Chandra}

Les $\cO$-blocs d'alg\`ebres de Hecke peuvent \^etre consid\'er\'es pour
d'autres choix de param\`etres. Dans cette partie, nous consid\'erons un
cas particuli\`erement important: celui des alg\`ebres de Hecke
associ\'ees aux ``groupes de Weyl'' relatifs
d'une $d$-s\'erie de Harish-Chandra de caract\`eres unipotents dans un
groupe r\'eductif fini. Nous \'etablis\-sons alors une compatibilit\'e
avec les familles de Lusztig.\par
On consid\`ere un groupe r\'eductif connexe sur une cl\^oture
alg\'ebrique d'un corps fini, muni d'un endomorphisme dont
une puissance est un Frobenius et on note $G$ le groupe des points fixes par
cet endomorphisme (un groupe fini ``de type de Lie'').
On note $\cE(G)$ l'ensemble de ses caract\`eres unipotents.
Soit $d>1$ un entier tel que le polyn\^ome cyclotomique
$\Phi_d(x)$ divise l'ordre g\'en\'erique de $G$. 
Suivant Brou\'e--Malle--Michel \cite[Thm.~3.2]{BMM0}, on consid\`ere la
partition de $\cE(G)$ en $d$-s\'eries de Harish-Chandra (cf. aussi
Brou\'e--Malle \cite{BM2})
\begin{equation} \label{HCdec}
  \cE(G)=\coprod_{L,\lambda} \cE(G,(L,\lambda))
\end{equation}
de $\cE(G)$, o\`u $L$ d\'ecrit l'ensemble des sous-groupes de Levi
$d$-ploy\'es de $G$ \`a conjugaison pr\`es, $\lambda$ d\'ecrit l'ensemble
des caract\`eres $d$-cuspidaux de $L$, et $\cE(G,(L,\lambda))$
consiste des constituants de $R_L^G(\lambda)$.

En outre, pour chaque paire cuspidale $(L,\lambda)$, on a une bijection
entre $\cE(G,(L,\lambda))$ et l'ensemble des caract\`eres irr\'eductibles
de l'alg\`ebre de Hecke $\cH(G,(L,\lambda))$ du
groupe de Weyl relatif $W_G(L,\lambda)$, avec certains param\`etres.
Ceci est d\'ecrit explicitement dans
\cite[Satz~3.14 et 6.10]{MaU} pour les groupes classiques et dans
\cite{BrMa} pour les groupes exceptionnels.
D'apr\`es Brou\'e--Malle--Michel \cite[Thm.~5.24]{BMM0}, la
partition~(\ref{HCdec}) est identique \`a la
partition de $\cE(G)$ en $\ell$-blocs, pour tout nombre premier $\ell$ divisant
$\Phi_d(q)$ mais ne divisant pas l'ordre du groupe de Weyl de $G$.\par
La partition en $\cO$-blocs (cf. chapitre~\ref{sec:theorie})
des caract\`eres irr\'eductibles de
chaque alg\`ebre de Hecke relative $\cH(G,(L,\lambda))$, fournit
une partition
$$\cE(G,(L,\lambda)) = \coprod_\cG \cG\,.$$
On a aussi une d\'ecomposition en familles de Lusztig de $\cE(G)$
\cite[Chap.~4]{LuB}. Pour les caract\`eres unipotents de la s\'erie
principale ($d=1$ et $L$ est un tore), elle co\"{\i}ncide avec la partition
pr\'ec\'edente en $\cO$-blocs \cite{Ro}.
\par
Ce r\'esultat se g\'en\'eralise au cas des paires cuspidales quelconques
et \`a tout $d$ en le th\'eor\`eme suivant:

\begin{thm} \label{fam=bloc} 
 Soit $G$ un groupe r\'eductif connexe sur un corps fini et
 $\cE(G,(L,\lambda))$ une $d$-s\'erie de Harish-Chandra de $\cE(G)$ d'alg\`ebre
 de Hecke relative $\cH(G,(L,\lambda))$. Supposons les
 conjectures~\ref{conj1}, \ref{conj2} et~\ref{conj3} vraies si le groupe de
 Weyl relative $W_G(L,\lambda)$ est de type exceptionnel non-r\'eel.
 Alors, les deux partitions
 suivantes de $\cE(G,(L,\lambda))$ co\"{\i}ncident:
 $$\coprod_{\cF} \Big(\cF\cap \cE(G,(L,\lambda)) \Big) = \coprod_\cG \cG\,,$$
 o\`u $\cF$ d\'ecrit les familles de Lusztig de $\cE(G)$ et
 $\cG$ d\'ecrit les $\cO$-blocs de $\Irr(\cH(G,(L,\lambda)))$.
 En d'autres termes, les $\cO$-blocs de
 $\Irr(\cH(G,(L,\lambda)))$ sont les intersections des familles de Lusztig avec
 $\cE(G,(L,\lambda))$. \par
 R\'eciproquement, la partition en familles de Lusztig de $\cE(G)$ est la
 plus fine partition telle que la propri\'et\'e pr\'ec\'edente soit
 vraie pour tout $d$ et toute paire $d$-cuspidale $(L,\lambda)$. 
\end{thm} 

\begin{proof}[Preuve]
Le reste de cette partie est consacr\'e \`a la preuve de ce th\'eor\`eme.
Notons tout d'abord qu'il suffit de consid\'erer le cas o\`u le
groupe alg\'ebrique est simple.

Consid\'erons pour commencer le cas des groupes classiques. Dans les
groupes de type $A_n$ et ${}^2A_n$, toutes les familles sont r\'eduites
\`a des singletons, tout comme les blocs des alg\`ebres de Hecke relatives
sur $\cO$.\par

Les caract\`eres unipotents des autres groupes classiques sont param\'etr\'es
par des symboles.
Un symbole $\Lambda=\{S,T\}$ est la donn\'ee de deux ensembles finis
$S$ et $T$ d'entiers positifs. On dit que deux symboles sont \'equivalents
si l'un peut \^etre obtenu \`a partir de l'autre par une suite d'op\'erations
de d\'ecalage
$$\{S,T\}\mapsto \{\{0\}\cup S+1,\{0\}\cup T+1\}\,.$$
Le d\'efaut d'un symbole est $\df(S,T)=\big||S|-|T|\big|$. Le rang de
$\Lambda=(S,T)$ est
$$\sum_{\lambda\in\Lambda} \lambda
  -\left\lfloor\frac{(|\Lambda|-1)^2}{4}\right\rfloor.$$
Pour $G$ de type $B_n$ ou $C_n$, alors $\cE(G)$ est param\'etr\'e par les
classes d'\'equivalence de symboles de rang $n$ et de d\'efaut impair. Deux
caract\`eres unipotents sont dans la m\^eme famille de Lusztig si et seulement
si les multi-ensembles r\'eunions de $S$ et $T$ co\"{\i}ncident.\par

Le param\'etrage des s\'eries de Harish-Chandra fourni par
\cite[Satz~3.14]{MaU} et la description des blocs pour les
alg\`ebres de Hecke donn\'ee par Brou\'e--Kim \cite[Thm.~3.16]{BrKi} montrent
que la premi\`ere partie du th\'eor\`eme est vraie.\par

Passons \`a la deuxi\`eme partie du th\'eor\`eme. Nous dirons qu'un 
sous-ensemble de $\cE(G)$ est une famille potentielle si elle est un
sous-ensemble minimal avec la propri\'et\'e que son intersection
avec toute $d$-s\'erie de Harish-Chandra est une union de
$\cO$-blocs pour l'alg\`ebre de Hecke relative correspondante.
Il s'agit donc de d\'emontrer que les familles potentielles sont les familles.
Notons pour commencer que la $d$-s\'erie de Harish-Chandra d'un symbole est
d\'etermin\'ee par le $d$-c{\oe}ur du symbole lorsque $d$ est impair et
par le $d/2$-coc{\oe}ur du symbole lorsque $d$ est pair \cite[Section 3A]{BMM0}.
En particulier, la $1$-s\'erie de Harish-Chandra d'un symbole est
d\'etermin\'ee par le d\'efaut du symbole.
On en d\'eduit que l'intersection d'une famille potentielle avec
une famille de Lusztig est une union de parties de la partition de cette
famille suivant les d\'efauts des symboles.
Il nous reste donc \`a d\'emontrer qu'\'etant donn\'es deux d\'efauts
$t_1$ et $t_2$ apparaissant dans une famille de Lusztig donn\'ee, alors
il existe un d\'efaut $d$ et deux symboles $\Lambda_1,\Lambda_2$
de d\'efauts $\df(\Lambda_i)=t_i$, tels que $\Lambda_1$ et $\Lambda_2$
sont dans la m\^eme $d$-s\'erie de Harish-Chandra, c'est-\`a-dire, ont
le m\^eme $d$-c{\oe}ur (respectivement $d/2$-coc{\oe}ur). \par

Soit donc $\cF$ une famille de Lusztig avec au moins deux \'el\'ements.
Le multi-ensemble des \'el\'ements de tout symbole de $\cF$ contient au moins
trois singletons.
Soit $\Lambda_1=\{S,T\}\in\cF$. On peut supposer qu'il existe
deux singletons $\lambda_1,\lambda_2$ dans $S$ avec
$\lambda_1<\lambda_2$.
Soit $\Lambda'$ le symbole obtenu \`a partir de $\Lambda_1$ en enlevant un
cocrochet de longueur $l:=\lambda_2-\lambda_1$ \`a la position
$\lambda_2$ de $S$ et soit $\Lambda_2$ le symbole obtenu \`a partir de
$\Lambda'$ en ajoutant un cocrochet de longueur $l$ \`a la position
$\lambda_1$ de $T$.
Alors, $\Lambda_2$ est obtenu \`a partir de $\Lambda_1$ en d\'epla\c{c}ant
$\lambda_1$ et $\lambda_2$ de $S$ \`a $T$, donc $\Lambda_2\in\cF$. Par
construction, $\Lambda_1$ et $\Lambda_2$ ont le m\^eme
$l$-coc{\oe}ur, donc appartiennent \`a la m\^eme 
$2l$-s\'erie de Harish-Chandra.
En outre, les d\'efauts de $\Lambda_1$ et $\Lambda_2$ diff\`erent de $4$.
Notons enfin que si $\cF$ contient des symboles de d\'efaut $t$, alors
tout entier impair de $\{1,\ldots, t\}$ appara\^{\i}t comme d\'efaut
d'un symbole de $\cF$. La construction pr\'ec\'edente permet donc d'atteindre
tous les d\'efauts de la famille.\par

Pour les groupes de type $D_n$ (respectivement 
${}^2D_n$), alors $\cE(G)$ est param\'etr\'e par des symboles de
d\'efaut congru \`a $0\pmod4$ (respectivement $2\pmod4$).
Des arguments similaires aux pr\'ec\'edents s'appliquent \`a cette situation.
\par
Passons maintenant aux types exceptionnels.
Les $d$-s\'eries de Harish-Chandra de groupes de Weyl relatifs cycliques
sont d\'ecrites dans Brou\'e--Malle \cite[Tab.~8.1 et ~8.3]{BrMa}. On d\'eduit
facilement de cette description la premi\`ere partie du th\`eor\`eme, dans
ce cas.

Rappelons que le cas de la $1$-s\'erie de Harish-Chandra principale est
d\'ej\`a connu par \cite{Ro}.

Dans la suite nous d\'eterminons les blocs des alg\`ebres de Hecke
associ\'ees aux groupes de r\'eflexions $W'$ de type exceptionnel pour les
param\`etres qui proviennent de cas o\`u $W'$ intervient comme groupe de Weyl
relatif d'une $d$-s\'erie de Harish-Chandra de caract\`eres unipotents d'un
groupe exceptionnel (voir table~\ref{tab:relatif}) et nous d\'emontrons la
premi\`ere partie du th\'eor\`eme~\ref{fam=bloc} pour ces cas.
Seul un groupe non cylique et non exceptionnel $W'$ appara\^{\i}t:
dans ce cas, les blocs sont d\'etermin\'es dans \cite{BrKi} et on d\'eduit
la premi\`ere partie du th\'eor\`eme.

Nous traitons aussi le cas du groupe de Coxeter $H_4$.

La deuxi\`eme partie du 
th\'eor\`eme~\ref{fam=bloc} pour les groupes $G$ exceptionnels
se r\'eduit alors \`a une v\'erification facile.
\end{proof}

\begin{table}[htbp] \caption{Groupes de Weyl relatifs non cycliques}
\label{tab:relatif}
\[\begin{array}{|r|rlllll|}
\hline
 G &d=\hfill 1& 3& 4& 5& 8& 12\cr
\hline
 {}^3D_4& G(6,6,2)& G_4& & & &\cr
 {}^2F_4& G(8,8,2)& -& G_{12}& & G_8& \cr
     F_4& G_{28},G(2,1,2)& G_5& G_8& & & \cr
     E_6& G_{35}& G_{25}& G_8& & & \cr
 {}^2E_6& G_{28}& G_5& G_8& & & \cr
     E_7& G_{36},G(2,1,3)& G_{26}& G_8,G(4,1,2)& & & \cr
     E_8& G_{37},G_{28},G(6,6,2)& G_{32}& G_{31}& G_{16}& G_9& G_{10}\cr
\hline
     H_4& G_{30}& G_{20}& G_{22}& & & \cr
\hline
\end{array}\]
\end{table}

\subsection{$W=G_4$}

Le groupe $W=G_4$ intervient comme groupe relatif dans les groupes r\'eductifs
de type $^3D_4$ (cf \cite[Folg.~5.6]{BrMa} pour les param\`etres).
Dans ce cas, l'alg\`ebre de Hecke a un seul bloc non trivial, de
caract\`eres $(\phi_{1,8},\phi_{2,1},\phi_{3,2})$. \par

\subsection{$W=G_5$}

Le groupe $G_5$ intervient comme groupe relatif dans les groupes r\'educ\-tifs
de type $F_4$, $E_6$ et $E_8$ (cf \cite[Table~5.9]{BrMa}). Nous indiquons
seulement les blocs non triviaux:

Dans $F_4$: $(\phi_{1,8}',\phi_{1,8}''',\phi_{2,5}'')$,
$(\phi_{1,12}',\phi_{1,12}'',\phi_{2,9})$,
$(\phi_{1,4}',\phi_{1,4}'',\phi_{1,16},\phi_{2,1},\phi_{2,7}',
 \phi_{2,7}'',\phi_{3,2},\ \ $
 $\phi_{3,4},\phi_{3,6})$.

Dans $E_6$: $(\phi_{1,4}',\phi_{1,16},\phi_{2,7}'')$,
$(\phi_{1,8}',\phi_{1,8}''',\phi_{2,5}'')$,
$(\phi_{1,12}',\phi_{2,3}',\phi_{3,2},\phi_{3,4},\phi_{3,6})$.

Dans $E_8$: $(\phi_{1,16},\phi_{2,1},\phi_{3,2},\phi_{3,4},\phi_{3,6})$.

\subsection{$W=G_8$}

Le groupe $G_8$ intervient comme groupe relatif dans les groupes r\'educ\-tifs
de type $^2F_4$, $F_4$, $E_6$, $E_7$ et $E_8$ (cf \cite[Table~5.12]{BrMa}).
Les blocs non triviaux sont:

Dans $^2F_4$: $(\phi_{2,1},\phi_{2,4})$, $(\phi_{2,10},\phi_{2,13})$,
$(\phi_{1,6},\phi_{1,12},\phi_{2,7}',\phi_{2,7}'',\phi_{3,4},\phi_{3,6},
  \phi_{4,3},\phi_{4,5})$.

Dans $F_4$: $(\phi_{2,1},\phi_{2,7}')$, $(\phi_{2,7}'',\phi_{2,13})$,
 $(\phi_{1,6},\phi_{1,18},\phi_{2,4},\phi_{2,10},\phi_{3,2},\phi_{3,6},
  \phi_{4,3},\phi_{4,5})$.

Dans $E_6$: $(\phi_{1,6},\phi_{2,4},\phi_{3,2})$,
$(\phi_{1,12},\phi_{2,10},\phi_{3,8})$,
$(\phi_{2,7}',\phi_{2,7}'',\phi_{4,3},\phi_{4,5})$.

Dans $E_7$: $(\phi_{2,1},\phi_{2,7}')$, $(\phi_{2,7}'',\phi_{2,13})$,
$(\phi_{3,2},\phi_{3,6})$, $(\phi_{4,3},\phi_{4,5})$,
$(\phi_{1,6},\phi_{1,18},\phi_{2,10})$.

Dans $E_8$: a) $(\phi_{2,1},\phi_{2,7}')$, $(\phi_{2,7}'',\phi_{2,13})$,
$(\phi_{3,2},\phi_{3,6})$, $(\phi_{4,3},\phi_{4,5})$,
$(\phi_{1,6},\phi_{1,18},\phi_{2,10})$;\par
b) $(\phi_{2,1},\phi_{2,7}'')$, $(\phi_{2,7}',\phi_{2,13})$,
$(\phi_{3,4},\phi_{3,8})$, $(\phi_{4,3},\phi_{4,5})$,
$(\phi_{1,0},\phi_{1,12},\phi_{2,4})$;\par
c) $(\phi_{2,1},\phi_{2,13},\phi_{4,3},\phi_{4,5})$;\par
d) $(\phi_{2,7}',\phi_{2,7}'',\phi_{4,3},\phi_{4,5})$.\par

\subsection{$W=G_9$}

Le groupe $G_9$ intervient comme groupe relatif dans les groupes r\'educ\-tifs
de type $E_8$ (cf \cite[Table~9]{MaD}). Les blocs non triviaux sont:

$(\phi_{2,7}',\phi_{2,11}'')$, $(\phi_{2,7}'',\phi_{2,10})$,
$(\phi_{1,12}',\phi_{1,12}'',\phi_{2,8},\phi_{2,11}')$, 
$(\phi_{1,6},\phi_{1,30},\phi_{2,13},\phi_{2,17})$,\par
$(\phi_{2,1},\phi_{2,4},\phi_{2,5},\phi_{2,14},\phi_{4,3},\phi_{4,5},
  \phi_{4,7},\phi_{4,9})$.

\subsection{$W=G_{10}$}

Le groupe $G_{10}$ intervient comme groupe relatif dans les groupes r\'eductifs
de type $E_8$ (voir \cite[Table~9]{MaD}). Les blocs non triviaux sont:

$(\phi_{1,6},\phi_{1,18})$, $(\phi_{1,22},\phi_{1,34})$, 
$(\phi_{2,1},\phi_{2,7}')$, $(\phi_{2,5},\phi_{2,11}')$,
$(\phi_{2,9},\phi_{2,15}')$, $(\phi_{2,7}'',\phi_{2,13})$,\par
$(\phi_{2,11}'',\phi_{2,17})$, $(\phi_{2,15}'',\phi_{2,21})$,
$(\phi_{4,3},\phi_{4,9})$, $(\phi_{4,7},\phi_{4,13})$,\par
$(\phi_{1,16},\phi_{2,10},\phi_{3,2},\phi_{3,6}'',\phi_{3,10}')$,
$(\phi_{1,12},\phi_{2,18},\phi_{3,6}',\phi_{3,10}'',\phi_{3,14})$,\par
$(\phi_{1,14},\phi_{1,26},\phi_{2,8},\phi_{3,12}',\phi_{3,12}'',\phi_{3,8}',
 \phi_{3,16},\phi_{3,8}'',\phi_{3,4},\phi_{4,11},\phi_{4,5})$,

Notons que pour montrer que 
$(\phi_{2,9},\phi_{2,15}')$ forme une famille, nous avons utilis\'e
\CHEVIE, pour calculer certaines valeurs de caract\`eres de l'alg\`ebre de
Hecke.

\subsection{$W=G_{12}$}

Le groupe $G_{12}$ intervient comme groupe relatif dans les groupes
r\'eductifs de type $^2F_4$ (cf \cite[Folg.~5.14]{BrMa}). Il y a un seul
bloc non trivial~:\par
$(\phi_{2,1},\phi_{2,4},\phi_{2,5},\phi_{4,3})$.

\subsection{$W=G_{16}$}

Le groupe $G_{16}$ intervient comme groupe relatif dans les groupes
r\'eductifs de type $E_8$ (cf \cite[Table~9]{MaD}). Les blocs non triviaux
sont:

$(\phi_{1,12},\phi_{2,7},\phi_{3,2})$, $(\phi_{1,36},\phi_{2,31},\phi_{3,26})$,
$(\phi_{3,10}',\phi_{3,10}'',\phi_{6,5})$,
$(\phi_{3,18}',\phi_{3,18}'',\phi_{6,13})$,\par
$(\phi_{2,13}',\phi_{2,13}'',\phi_{4,3},\phi_{4,8})$,
$(\phi_{2,25}',\phi_{2,25}'',\phi_{4,15},\phi_{4,20})$,
$(\phi_{4,6},\phi_{4,11})$, $(\phi_{4,12},\phi_{4,17})$,\par
$(\phi_{1,24},\phi_{2,19}',\phi_{2,19}'',\phi_{3,14}',\phi_{3,14}'',\phi_{4,9},
  \phi_{4,14},\phi_{5,4},\phi_{5,8},\phi_{5,10},\phi_{5,12},\phi_{5,16},
  \phi_{6,9})$.

\subsection{$W=G_{20}$}

Le groupe $G_{20}$ intervient comme groupe relatif pour le spets de type
$H_4$ (cf \cite[Table~9]{MaD}). Les blocs non triviaux sont:

$(\phi_{2,1},\phi_{2,7})$, $(\phi_{2,21},\phi_{2,27})$,
$(\phi_{4,3},\phi_{4,6})$, $(\phi_{4,13},\phi_{4,16})$,\par
$(\phi_{1,20},\phi_{2,11},\phi_{2,17},\phi_{3,2},\phi_{3,6},\phi_{3,10}',
  \phi_{3,10}'',\phi_{3,12},\phi_{3,14},\phi_{4,8},\phi_{4,11},\phi_{5,8},
  \phi_{6,5},\phi_{6,7},\phi_{6,9})$.

\subsection{$W=G_{22}$}

Le groupe $G_{22}$ intervient comme groupe relatif pour le spets de type
$H_4$ (cf \cite[Table~9]{MaD}). Les blocs non triviaux sont:

$(\phi_{3,2},\phi_{3,6})$, $(\phi_{3,12},\phi_{3,16})$,
$(\phi_{2,1},\phi_{2,7},\phi_{2,11},\phi_{2,13},\phi_{4,3},\phi_{4,6},
  \phi_{4,8},\phi_{4,9},\phi_{6,5},\phi_{6,7})$.

\subsection{$W=G_{25}$}

Le groupe $G_{25}$ intervient comme groupe relatif dans les groupes
r\'eductifs de type $E_6$ (cf \cite[Folg.~5.16]{BrMa}). Les blocs
non triviaux sont:

$(\phi_{3,5}',\phi_{3,5}'',\phi_{6,2})$,
$(\phi_{3,13}',\phi_{3,13}'',\phi_{6,10})$,
$(\phi_{1,12},\phi_{2,9},\phi_{3,6},\phi_{8,6},\phi_{9,5},\phi_{9,7})$.

\subsection{$W=G_{26}$}

Le groupe $G_{26}$ intervient comme groupe relatif dans les groupes
r\'eductifs de type $E_7$ (cf. \cite[Prop.~7.1]{MaG}). Les blocs non triviaux
sont:

$(\phi_{1,12},\phi_{2,9},\phi_{3,6})$,
$(\phi_{1,21},\phi_{2,18},\phi_{3,15})$,
$(\phi_{3,5}',\phi_{3,5}'',\phi_{6,2})$,
$(\phi_{3,16}',\phi_{3,16}'',\phi_{6,13})$,\par
$(\phi_{3,8}',\phi_{3,8}'',\phi_{6,5})$,
$(\phi_{3,13}',\phi_{3,13}'',\phi_{6,10})$, $(\phi_{8,6}'',\phi_{8,9}')$,
$(\phi_{1,9},\phi_{2,15},\phi_{8,3},\phi_{8,6}',\phi_{9,5},\phi_{9,7})$,\par
$(\phi_{1,24},\phi_{2,12},\phi_{8,9}'',\phi_{8,12},\phi_{9,8},\phi_{9,10})$.

\subsection{$W=G_{28}$}

Le groupe $G_{28}=W(F_4)$ intervient comme groupe relatif dans les groupes
r\'eductifs de type $E_6$ et $E_8$. Les blocs non triviaux sont:

Dans $E_6$: $(\phi_{1,12}',\phi_{2,4}'',\phi_{9,2},\phi_{8,3}')$,
$(\phi_{1,12}'',\phi_{2,16}',\phi_{9,10},\phi_{8,9}'')$,
$(\phi_{4,8},\phi_{6,6}',\phi_{6,6}'',\phi_{12,4},\phi_{16,5})$.\par
Dans $E_8$: 
$(\phi_{4,8},\phi_{6,6}',\phi_{6,6}'',\phi_{12,4},\phi_{16,5})$.

\subsection{$W=G_{31}$}
Le groupe $G_{31}$ intervient comme groupe relatif dans les groupes
r\'eductifs de type $E_8$ (cf. \cite[Prop.~7.1]{MaG}). Les blocs non triviaux
sont:

$(\phi_{4,1},\phi_{4,7})$, $(\phi_{10,2},\phi_{10,6})$,
$(\phi_{20,3},\phi_{20,5})$, $(\phi_{20,7},\phi_{20,13}')$,
$(\phi_{36,5},\phi_{36,7})$, $(\phi_{40,7},\phi_{40,9})$,\par
$(\phi_{40,13},\phi_{40,15})$, $(\phi_{36,15},\phi_{36,17})$,
$(\phi_{20,13}'',\phi_{20,19})$, $(\phi_{20,21},\phi_{20,23})$,
$(\phi_{10,26},\phi_{10,30})$,\par $(\phi_{4,31},\phi_{4,37})$,
$(\phi_{10,12},\phi_{15,8}',\phi_{15,8}'',\phi_{30,4},\phi_{40,6})$,
$(\phi_{10,24},\phi_{15,20}',\phi_{15,20}'',\phi_{30,16},\phi_{40,18})$,\par
$(\phi_{6,18},\phi_{6,14},\phi_{16,16},\phi_{24,14},\phi_{24,6},\phi_{30,10}',
  \phi_{30,10}'',\phi_{36,10},\phi_{40,10},\phi_{40,14},\phi_{64,9},
  \phi_{64,11})$.

\subsection{$W=G_{32}$}

Le groupe $G_{32}$ intervient comme groupe relatif dans les groupes
r\'eductifs de type $E_8$ (cf. \cite[Prop.~7.1]{MaG}). Les blocs non triviaux
sont:

$(\phi_{64,8},\phi_{64,11})$, $(\phi_{64,18},\phi_{64,21})$,\par
$(\phi_{4,11},\phi_{6,8},\phi_{10,2})$,
$(\phi_{4,51},\phi_{6,48},\phi_{10,42})$,
$(\phi_{20,13},\phi_{40,10},\phi_{60,7})$,
$(\phi_{20,25},\phi_{40,22},\phi_{60,19})$,\par
$(\phi_{5,12},\phi_{15,6},\phi_{20,3})$,
$(\phi_{5,44},\phi_{15,38},\phi_{20,35})$,
$(\phi_{10,10},\phi_{20,7},\phi_{30,4})$,
$(\phi_{10,34},\phi_{20,31},\phi_{30,28})\!$,\par
$(\phi_{20,17},\phi_{40,14},\phi_{60,11}'')$,
$(\phi_{20,21},\phi_{40,18},\phi_{60,15}'')$,
$(\phi_{20,12},\phi_{60,12},\phi_{80,9})$,\par
$(\phi_{20,20},\phi_{60,20},\phi_{80,17})$,
$(\phi_{30,12}',\phi_{30,12}'',\phi_{60,9})$,
$(\phi_{30,20}',\phi_{30,20}'',\phi_{60,17})$,\par
$(\phi_{4,21},\phi_{20,9}',\phi_{20,9}'',\phi_{24,6},\phi_{36,7},\phi_{36,5})$,
$(\phi_{4,41},\phi_{20,29}',\phi_{20,29}'',\phi_{24,26},\phi_{36,27},\phi_{36,25})$,\par
$(\phi_{5,20},\phi_{10,14},\phi_{15,8},\phi_{40,8},\phi_{45,6},\phi_{45,10})$,
$(\phi_{5,36},\phi_{10,30},\phi_{15,24},\phi_{40,24},\phi_{45,22},\phi_{45,26})$,\par
$(\phi_{1,40},\phi_{6,28},\phi_{15,16},\phi_{15,22},\phi_{20,19},\phi_{20,16},
\phi_{24,16},\phi_{36,17},\phi_{36,15},\phi_{45,14},\phi_{45,18},\phi_{60,13},$\par
$\phi_{60,16},\phi_{64,16},\phi_{64,13},\phi_{80,13},\phi_{81,12},\phi_{81,14},
\phi_{81,10})$.


\end{document}